\documentclass{article}

\usepackage{arxiv}

\usepackage[utf8]{inputenc} % allow utf-8 input
\usepackage[T1]{fontenc}    % use 8-bit T1 fonts
\usepackage{hyperref}       % hyperlinks
\usepackage{url}            % simple URL typesetting
\usepackage{booktabs}       % professional-quality tables
\usepackage{amsfonts}       % blackboard math symbols
\usepackage{nicefrac}       % compact symbols for 1/2, etc.
\usepackage{microtype}      % microtypography
\usepackage{lipsum}		% Can be removed after putting your text content
\usepackage{graphicx}
\usepackage{natbib}
\usepackage{doi}

\usepackage{amssymb}
\usepackage{latexsym}

\usepackage{graphicx}
\usepackage{amsmath}
\usepackage{bm}
\usepackage{caption}
\usepackage{subcaption}
\usepackage[version=4]{mhchem}
\usepackage{siunitx}
\usepackage{longtable,tabularx}
\newcommand{\sgn}{\mathop{\mathrm{sgn}}}

\title{Indirect Optimization of Multi-Phase Trajectories Involving Arbitrary Discrete Logic}

\author{Harish Saranathan \\
        Independent Researcher \\
	Millis, MA 02054, USA \\
	\texttt{hsaranat@alumni.purdue.edu} \\
}

\date{}

\hypersetup{
pdftitle={Indirect Optimization of Multi-Phase Trajectories Involving Arbitrary Discrete Logic},
pdfsubject={math.OC},
pdfauthor={Harish Saranathan},
pdfkeywords={trajectory optimization, indirect methods, entry, descent, and landing, hybrid systems, mixed integer programming},
}

\begin{document}
\maketitle

\begin{abstract}
Multi-phase trajectories of aerospace vehicle systems involve multiple flight segments whose transitions may be triggered by boolean logic in continuous state variables, control and time. When the boolean logic is represented using only states and/or time, such systems are termed autonomously switched hybrid systems. The relaxed autonomously switched hybrid system approach (RASHS) was previously introduced to simplify the trajectory optimization process of such systems in the indirect framework when the boolean logic is solely represented using AND operations. This investigation enables cases involving arbitrary discrete logic. The new approach is termed the Generalized Relaxed Autonomously Switched Hybrid System (GRASHS) approach. Similar to the RASHS approach, the outcome of the GRASHS approach is the transformation of the necessary conditions of optimality from a multi-point boundary value problem to a two-point boundary value problem, which is simpler to handle. This is accomplished by converting the arbitrary boolean logic to the disjunctive normal form and applying smoothing using sigmoid and hyperbolic tangent functions. The GRASHS approach is demonstrated by optimizing a Mars entry, descent, and landing trajectory, where the parachute descent segment is active when the velocity is below the parachute deployment velocity or the altitude is below the parachute deployment altitude, and the altitude is above the powered descent initiation altitude. This set of conditions represents a combination of AND and OR logic. The previously introduced RASHS approach is not designed to handle such problems. The proposed GRASHS approach aims to fill this gap.
\end{abstract}

% keywords can be removed

{\renewcommand\arraystretch{1.0}
\noindent\begin{longtable*}{@{}l @{\quad=\quad} l@{}}
CSC & composite smooth control \\
EDL & entry, descent, and landing \\
DAE & differential algebraic equation \\
DNF & disjunctive normal form \\
DOF & degrees-of-freedom \\
ENU & east-north-up \\
GRASHS & generalized relaxed autonomously switched hybrid system \\
HTS & hyperbolic tangent smoothing \\
MPBVP & multi-point boundary value problem \\
MSL & mars science laboratory \\
PCPF & planet-centered planet-fixed \\
PDI & powered descent initiation \\
RASHS & relaxed autonomously switched hybrid system \\
SCP & sequential convex programming \\
SOP & sum of products \\
STC & state-triggered constraint \\
TPBVP & two-point boundary value problem \\
$\texttt{A}, \texttt{B}, \texttt{C}, \texttt{D}$ & logical variables \\
$B$ & logical expression \\
$C_D$ & drag coefficient \\
$C_L$ & lift coefficient \\
$D$ & drag, N \\
$\bm{\hat{\textbf{e}}}_E$, $\bm{\hat{\textbf{e}}}_N$, $\bm{\hat{\textbf{e}}}_Z$ & unit vectors defining local East-North-Up frame \\
$\textbf{F}$ & force vector, N \\
$\textbf{f}$  & equations of motion \\
$g$ & condition on continuous state vector and time \\
$\textbf{H}$ & Hamiltonian \\
$H$ & atmospheric scale height, m \\
$h$ & altitude, h \\
$\hbar$ & scaled altitude, nd \\
$I_{sp}$ & specific impulse, s \\
$J$ & cost functional \\
$K_1$, $K_2$, $K_3$ & weights \\
$k$ & empirical constant for heat-rate calculation, J/s\textsuperscript{2}.kg\textsuperscript{-1/2}.m\textsuperscript{-3} \\
$L$ & lift, N \\
$\mathcal{L}$ & Lagrangian \\
$M$ & scaled mass, nd \\
$m$ & number of trajectory segments or mass, kg \\
$\dot{m}$ & mass flow rate, kg/s \\
$n$ & number of conditions \\
$p$ & number of minterms \\
$Q$ & stagnation-point heat-load, J/m\textsuperscript{2} \\
$q$ & number of variables in minterm \\
$\dot{q}$ & stagnation-point heat-rate, W/m\textsuperscript{2} \\
$R$ & radius of planet, m \\
$R_N$ & nose radius, m \\
$\textbf{r}$ & inertial position vector, m \\
$S$ & reference area, m\textsuperscript{2} \\
$s$ & slope of sigmoid function \\
$T$ & thrust, N \\
$t$ & time, s \\
$\textbf{U}$ & control vector \\
$u$ & unit step function \\
$V$ & scaled atmospheric-relative velocity, nd \\
$\textbf{v}$ & atmospheric-relative velocity vector, m/s \\
$v$ & atmospheric-relative velocity, m/s \\
$\textbf{X}$ & continuous state vector \\
$\bm{\hat{\textbf{X}}}$, $\bm{\hat{\textbf{Y}}}$, $\bm{\hat{\textbf{Z}}}$ & unit vectors defining inertial frame \\
$\bm{\hat{\textbf{X}}_G}$, $\bm{\hat{\textbf{Y}}_G}$, $\bm{\hat{\textbf{Z}}_G}$ & unit vectors defining planet-centered planet-fixed frame \\
$\bm{\hat{\textbf{x}}}_B$, $\bm{\hat{\textbf{y}}}_B$, $\bm{\hat{\textbf{z}}}_B$ & unit vectors defining body frame \\
$\bm{\hat{\textbf{x}}}_W$, $\bm{\hat{\textbf{y}}}_W$, $\bm{\hat{\textbf{z}}}_W$ & unit vectors defining wind frame \\
$\texttt{z}$ & generic function \\
$\alpha$ & angle-of-attack, rad \\
$\gamma$ & atmospheric-relative flight path angle \\
$\zeta$ & slope parameter of hyperbolic tangent function \\
$\theta$ & longitude, rad \\
$\bm{\lambda}$ & co-state vector \\
$\mu$ & standard gravitational parameter, m\textsuperscript{3}/s\textsuperscript{2} \\
$\xi$ & switching function \\
$\bm{\Pi}$ & Lagrange multiplier vector adjoining interior-point boundary conditions \\
$\rho$ & atmospheric density, kg/m\textsuperscript{3} \\
$\sigma$ & bank angle, rad \\
$\tau$ & proxy variable for time used in integration \\
\textbf{$\Phi$} & terminal cost \\
$\phi$ & latitude, rad \\
$\bm{\Psi}$ & boundary condition vector \\
$\psi$ & heading angle, rad \\
$\bm{\omega}$ & angular velocity of planet, rad/s \\

\multicolumn{2}{@{}l}{Subscripts}\\
$F$ & fuel \\
$f$ & final time \\
$i$ & first index of conditions in trajectory segment\\
$j$ & second index of conditions in trajectory segment\\
$k$ & trajectory segment index\\
$max$ & maximum \\
$o$ & surface \\
$P$ & parachute deployment \\
$PATH$ & path \\
$PDI$ & powered descent initiation \\
$p$ & product term in DNF \\
$0$ & initial time \\
\multicolumn{2}{@{}l}{Superscripts}\\
$i$ & inertial
\end{longtable*}}

\section{Introduction}
Aerospace vehicles employed in missions such as entry, descent, and landing (EDL) fly multi-phase trajectories consisting of multiple flight segments. An EDL trajectory typically consists of a hypersonic, parachute descent and powered descent segment. Vehicles flying such trajectories are categorized as hybrid systems (\cite{literature1}) because their motion is described by continuous states (such as position and velocity) and a discrete state that governs which flight segment is active (also referred to as the mode of operation) at any given time. The system is autonomously switched when the discrete state is solely dependent on the continuous state variables and/or time, such as in an EDL system, where each flight segment transition is typically triggered by conditions on velocity and altitude. Conversely, in a non-autonomously switched hybrid system, the mode of operation is explicitly controlled using a combination of continuous and discrete control inputs.

In a multi-phase flight trajectory, each segment is subject to a set of equations of motion that are continuous and differentiable in that segment. These equations typically have discontinuities at the transition points of the flight segments because of discrete changes in physical characteristics of the vehicle such as mass and aerodynamic coefficients. Therefore, for the overall trajectory, the equations of motion are piecewise continuous.

The optimization of such trajectories the minimization (or maximization) of a cost functional that consists of a path cost, which can be different for each segment. The path cost of a given segment beginning at time $t_1$ and ending at time $t_2$ is typically represented as:

\begin{equation}
    J_{PATH} = \int\limits_{t_{1}}^{t_2} \mathcal{L} d\tau
\end{equation}

\noindent where $\mathcal{L}$ is the Lagrangian. For an EDL trajectory, the path cost during hypersonic segment may be the stagnation-point heat-load, while that for powered descent segment may be the total propellant consumed. Consequently, the Lagrangian for the hypersonic phase would be the stagnation-point heat-rate, while that for the powered descent segment would be the thrust. Therefore, the Lagrangian for the overall multi-phase trajectory is also piecewise continuous.

Typical trajectory optimization techniques for non autonomously switched hybrid systems, which is beyond the scope of this investigation, are summarized by \cite{literature2}, \cite{literature3}, \cite{literature4} and \cite{literature5}. Trajectories of autonomously switched hybrid systems, which is the focus of this investigation, have traditionally been optimized using direct methods (\cite{literature6, literature7, literature8}). Conversely, indirect methods have historically not been extensively adopted when dealing with such systems because of associated challenges. When employing indirect methods, the optimization problem is transcribed into a multi-point boundary value problem (MPBVP) that represents the necessary conditions of optimality in a system of differential algebraic equations (DAE) (\cite{literature9, literature10}). Solving this system of equations poses challenges because existing numerical algorithms require an initial guess for each segment that must already be close to the actual solution. Additionally, as the number of segments increases, the number of interior-point boundary conditions at flight segment transition points that must be enforced also increases. A common solution strategy in the indirect framework involved a mixed approach (\cite{literature11, literature12}), wherein the switching instants and the continuous states at these instants were parameterized and optimized using techniques such as gradient methods, and the trajectories between these switching instants were optimized using indirect methods. This approach becomes exponentially complex with the number of flight segments because it requires the evaluation of every possible sequence of flight segments.

\cite{rashs, rashs2, rashs3} and \cite{continuation4} introduced the relaxed autonomously switched hybrid system (RASHS) approach to mitigate the challenges associated with solving the resultant MPBVP in the indirect framework. The approach relaxes the MPBVP using sigmoid functions to approximate the piecewise equations of motion and Lagrangian into continuous equations. This reduces the MPBVP to a two-point boundary value problem (TPBVP), wherein only the end-point boundary conditions must be explicitly enforced. The intermediate boundary conditions at the flight segment transition points are implicitly satisfied by the new continuous equations of motion and cost functional.

Motivated in part by RASHS, \cite{csc1, csc2} recently developed the Composite Smooth Control (CSC) framework to handle non-autonomously switched hybrid systems in the indirect framework. In this approach, the equations of motion and cost functional are smoothed using sigmoid functions and the discrete control inputs are smoothed using the hyperbolic tangent smoothing (HTS) technique (\cite{hts}).

However, the RASHS and CSC techniques can only be applied to problems where a given flight segment is active when every condition associated with that segment is satisfied. That is, the conditions associated with a given segment must be solely represented using AND logic. The limitation of this formulation for autonomously switched hybrid systems is that the switching conditions must be appropriately chosen to guarantee that every flight segment is activated during the mission. For example, consider the aforementioned EDL mission consisting of the three flight segments, with conditions governing the active segment described in Table \ref{tab:EDLflightSegments}. The quantities $h$ and $v$ represent altitude and atmospheric-relative velocity respectively. The subscripts $P$ and $PDI$ denote that the specified quantity corresponds to parachute deployment and powered descent initiation (PDI) respectively. Note that the conditions in Table \ref{tab:EDLflightSegments} only involve AND logic. If the chosen value of parachute deployment velocity, $v_P$, is low, the vehicle may never decelerate to this velocity prior to descending to PDI altitude, $h_{PDI}$. As a result, the parachute descent segment will be skipped altogether.

\begin{table}[!h]
    \centering
    \caption{\label{tab:EDLflightSegments}EDL flight segments and associated conditions.}
    \begin{tabular}{lccc}
        \hline
        Flight Segment & Conditions
        \\\hline
        Hypersonic to low-supersonic & $v \geq v_P$ \\
        Parachute descent & $v < v_P$ AND $h \geq h_{PDI}$ \\
        Powered descent & $h < h_{PDI}$ \\
        \hline
    \end{tabular}
\end{table}

Choosing a reasonable value for the parachute deployment velocity to guarantee the activation of parachute descent segment requires prior analysis of the flight dynamics. This approach can become time-consuming as the number of flight segments and associated conditions increase. Therefore, it may become necessary to implement a guard condition using OR logic. For instance, the conditions can be formulated such that the parachute descent segment is activated when either the velocity reduces below the parachute deployment velocity OR the altitude reduces below the parachute deployment altitude. The additional condition on altitude serves as a contingency in the event the vehicle does not decelerate to the parachute deployment velocity prior to descending to the PDI altitude. Therefore, the conditions associated with parachute descent segment can be represented as $\left(v < v_P \textrm{ OR } h < h_P\right) \textrm{ AND } h \geq h_{PDI}$. However, RASHS and CSC formulations are not set up to handle such conditions because they entail OR logic. Instead, these formulations require assumptions about whether the vehicle first decelerates to $v_P$ or descends to $h_P$ to eliminate the OR logic. Consequently, this involves trial and error, wherein the altitude trigger may first be ignored. If the parachute descent segment never gets activated, it can be inferred that the vehicle descends to $h_P$ before decelerating to $v_P$. Therefore, the problem must be solved again by ignoring the velocity trigger as opposed to the altitude trigger. In essence, in the worst case scenario, this problem will have to be solved twice. This count will exponentially increase with the number of flight segments and the associated conditions involving OR logic. To avoid solving the problem multiple times, it would be beneficial to develop a method that can simultaneously handle a combination of AND and OR logic.

In the direct trajectory optimization arena, \cite{malyuta} developed a fast homotopy approach to handle conditions represented using OR logic and applied it to spacecraft rendezvous trajectory optimization. In this approach, the discrete OR logic is smoothed using a multinomial logit function (\cite{logit}) and embedded into a sequential convex programming (SCP) framework using continuous embedding (\cite{literature3}). There have also been recent advances in SCP employing continuous state-triggered constraints (STCs) to handle discrete logic (\cite{stcscp}). However, it must be noted that these methods are developed for direct trajectory optimization through SCP.

Therefore, a gap continues to exist in the indirect trajectory optimization arena because of the lack of a framework that can handle arbitrary discrete logic. Motivated by the RASHS framework, the work presented in this investigation, termed the Generalized Relaxed Autonomously Switched Hybrid System (GRASHS) approach, serves to fill this gap specifically for autonomously switched hybrid systems. In this approach, the arbitrary discrete logic is transformed to the disjunctive normal form (DNF) (\cite{dnf}). A DNF represents any arbitrary boolean logic solely using AND, OR and NOT operations. The GRASHS approach then represents the AND operation as a product of boolean values and the OR operation as the signum of the sum of boolean values. The NOT operation is trivially handled by replacing the operation with appropriate predicates. The AND logic is smoothed using sigmoid functions and the OR logic is smoothed using the hyperbolic tangent function. The resultant equations of motion and Lagrangian are smooth for the entire trajectory, thereby reducing the MPBVP representing the necessary conditions of optimality to a TPBVP. This framework is demonstrated using a Mars EDL example, where the conditions activating the parachute descent segment are represented using a combination of AND and OR logic. The problem is solved with no apriori knowledge about which condition constituting the OR logic triggers the parachute descent segment. The results are compared against the solutions of the RASHS formulation and the original MPBVP. These comparison solutions are calculated by eliminating the OR logic using additional knowledge gained from the GRASHS solution.

\section{The Generalized Relaxed Autonomously Switched Hybrid System (GRASHS) Approach}
\label{sec:grashs}
Consider a multi-phase trajectory consisting of $m$ flight segments, where segment $k$ is active if a boolean expression $B_k$ represented using an arbitrary combination of conditions, each bearing the form $g\left(\textbf{X},t\right) < 0$, is \texttt{true}. Here on, the conditional expressions $g < 0$ will be termed predicates. Because the system is autonomously switched, the predicates are solely dependent on the continuous state vector $\textbf{X}$ and time $t$. Let the segment $k$ be governed by the equations of motion $\textbf{f}_k$. Let the Lagrangian of the path cost to be minimized in that segment be $\mathcal{L}_k$. Therefore:

\begin{equation}
    \label{eqn:rashs_singleExp}
    \begin{gathered}
        \dot{\textbf{X}} = \sum_{k=1}^{m} \xi_k \textbf{f}_k \\
        \mathcal{L} = \sum_{k=1}^{m} \xi_k \mathcal{L}_k
    \end{gathered}
\end{equation}

\noindent where $\xi_k$ is a switching function (the discrete state) defined as:

\begin{equation}
    \label{eqn:rashs_switchingFunction}
    \xi_k =
    \begin{cases}
        1 \text{ when } B_k = \texttt{true} \\
        0 \text{ otherwise}
    \end{cases}
\end{equation}

The goal of the GRASHS approach is to convert Eq. \eqref{eqn:rashs_singleExp} into a smooth approximation when $B_k$ is represented using an arbitrary combination of boolean operations on the predicates. As a first step, $B_k$ must be converted into the disjunctive normal form (DNF) (\cite{dnf}), which is one of the standard representations of any boolean expression. The DNF consists of a sum (OR) of one or more product terms (AND), also known as the minterms, and inversion (NOT operation) of some boolean variables. Therefore, the DNF is sometimes also referred to as the sum of products (SOP). The following is an example of a DNF expression consisting of logical variables $\texttt{A}$, $\texttt{B}$, $\texttt{C}$ and $\texttt{D}$.

\begin{equation}
    \label{eqn:sopExample}
    \overline{\texttt{A}} \cdot \texttt{B} \cdot \texttt{C} \cdot \overline{\texttt{D}}
    +
    \texttt{A} \cdot \overline{\texttt{B}} \cdot \texttt{D}
    +
    \texttt{C} \cdot \overline{\texttt{D}}
    +
    \texttt{A} \cdot \overline{\texttt{D}}
\end{equation}

In Eq. \eqref{eqn:sopExample}, the dot ($\cdot$) represents AND, the plus ($+$) represents OR and the bar ($\overline{\phantom{x}}$) represents inverter (NOT). The minterms in this equation are ($\overline{\texttt{A}} \cdot \texttt{B} \cdot \texttt{C} \cdot \overline{\texttt{D}}$), ($\texttt{A} \cdot \overline{\texttt{B}} \cdot \texttt{D}$), ($\texttt{C} \cdot \overline{\texttt{D}}$) and ($\texttt{A} \cdot \overline{\texttt{D}}$). In the autonomously switched hybrid system trajectory optimization problem, the logical variables $\texttt{A}$, $\texttt{B}$, $\texttt{C}$ and $\texttt{D}$ represent the predicates $g_{i,j,k} < 0$ that constitute $B_k$. Note that subscript $i$, $j$, and $k$ have been added. The subscript $i$ represents the index of the minterm, $j$ represents the index of the predicate inside a given minterm, and $k$ represents the index of the flight segment. When the variables are inverted (inversion represents a NOT operation), the predicates become $g_{i,j,k} >= 0$. However, this expression can simply be replaced by $-g_{i,j,k} < 0$. This way, the GRASHS approach can be trivially extended to handle NOT operations. Conversion of an arbitrary boolean expression into DNF is explained by \cite{dnf}.

In the EDL mission such as that described in Table \ref{tab:EDLflightSegments}, the parachute descent segment is the second segment ($k=2$). Suppose this segment is active when velocity $v$ is less than the parachute deployment velocity $v_P$ OR the altitude $h$ is less than the parachute deployment altitude $h_P$, AND the altitude is greater than or equal to the powered descent initiation altitude $h_{PDI}$, the corresponding predicates are as follows:

\begin{equation}
    \label{eqn:grashs_EDLExample1}
    \begin{gathered}
        g_{1,1,2} < 0 \Longrightarrow v - v_P < 0 \\
        g_{1,2,2} < 0 \Longrightarrow h_{PDI} - h < 0 \\
        g_{2,1,2} < 0 \Longrightarrow h - h_P < 0 \\
        g_{2,2,2} < 0 \Longrightarrow h_{PDI} - h < 0
    \end{gathered}
\end{equation}

The logical expression $B_k$ ($k=2$) can then be represented in DNF as follows:

\begin{equation}
    \label{eqn:grashs_EDLExample2}
    \left(g_{1,1,2} < 0 \right) \cdot \left(g_{1,2,2} < 0\right) + \left(g_{2,1,2} < 0\right) \cdot \left(g_{2,2,2} < 0\right)
\end{equation}

Note that although $g_{1,2,2}$ and $g_{2,2,2}$ are essentially the same, they are listed as separate predicates because they belong to different minterms. This makes it easier to generalize Eq. \eqref{eqn:grashs_EDLExample2}. As a generalization of Eq. \eqref{eqn:grashs_EDLExample2}, suppose for flight segment $k$, there are $p_k$ minterms and a given minterm consists of $q_{p_k}$ predicates, then $B_k$ can be represented as a DNF as follows:

\begin{equation}
    \label{eqn:grashs_Bk}
    B_k\left(\textbf{X},t\right) = 
    \overset{p_k}{\underset{i=1}{\lor}} \left(\overset{q_{p_k}}{\underset{j=1} {\land}}\left(g_{i,j,k}\left(\textbf{X},t\right)<0\right)\right)
\end{equation}

\noindent where $\land$ and $\lor$ represent AND and OR operations respectively. The AND expression $\overset{q_{p_k}}{\underset{j=1} {\land}}\left(g_{i,j,k}\left(\textbf{X},t\right)<0\right)$ represents the $i$\textsuperscript{th} minterm. Accordingly, Eq. \eqref{eqn:rashs_switchingFunction} becomes:

\begin{equation}
    \label{eqn:grashs_switchingFunction}
    \xi_k =
    \begin{cases}
        1 \text{ when } \overset{p_k}{\underset{i=1}{\lor}} \left(\overset{q_{p_k}}{\underset{j=1} {\land}}\left(g_{i,j,k}\left(\textbf{X},t\right)<0\right)\right) = \texttt{true} \\
        0 \text{ otherwise}
    \end{cases}
\end{equation}

Noting that the predicates $g_{i,j,k} < 0$ evaluate to either \texttt{true} or \texttt{false}, each predicate can be represented using a horizontally flipped unit step function of the form $\left[1 - u\left(g_{i,j,k}\left(\textbf{X},t\right)\right)\right]$. This expression evaluates to $1$ when the predicate is \texttt{true} and evaluates to $0$ otherwise. Because the AND expression representing a given minterm evaluates to \texttt{true} only if every predicate in it is \texttt{true}, the minterm can be expressed as a product of the horizontally flipped unit step functions as follows:

\begin{equation}
    \label{eqn:GRASHSMintermStep}
    \overset{q_{p_k}}{\underset{j=1} {\land}}\left(g_{i,j,k}\left(\textbf{X},t\right)<0\right)
    \Longrightarrow
    \prod_{j=1}^{q_{p_k}} \left[1 - u\left(g_{i,j,k}\left(\textbf{X},t\right)\right)\right]
\end{equation}

The product equals $1$ only when every horizontally flipped unit step function evaluates to $1$, implying that every predicate associated with the minterm must evaluate to \texttt{true}. Otherwise, the product evaluates to $0$. To implement the OR logic on the minterms represented by Eq. \eqref{eqn:grashs_switchingFunction}, the following must be noted:

\begin{itemize}
    \item Each minterm represented by Eq. \eqref{eqn:GRASHSMintermStep} will evaluate to either $0$ (\texttt{false}) or $1$ (\texttt{true}).
    \item The overall OR logic involving all minterms must evaluate to $1$ when at least one of the minterms evaluates to $1$.
    \item The overall OR logic involving all minterms must evaluate to $0$ only when every minterm evaluates to $0$.
\end{itemize}

Also, note that:
\begin{itemize}
    \item The summation of the minterms will evaluate to $0$ if every minterm evaluates to $0$.
    \item The summation of the minterms will evaluate to a positive integer if at least one minterm evaluates to $1$.
\end{itemize}

Therefore, the OR logic can be represented as a summation of the minterms of the form in Eq. \eqref{eqn:GRASHSMintermStep} and saturating the summation at $1$. The saturation can be accomplished by applying the signum function on the summation. The signum function evaluates to $1$ when its input is positive, $0$ when the input is $0$, and $-1$ when the input is negative. The negative input is not applicable because the summation of minterms will always be nonnegative. Consequently, Eq. \eqref{eqn:sopExample} becomes:

\begin{equation}
    \label{eqn:grashs_unitstep}
    \begin{gathered}
        \dot{\textbf{X}} = \sum_{k=1}^{m} \left( \sgn\left(\sum_{i=1}^{p_k} \left(\prod_{j=1}^{q_{p_k}} \left[1 - u\left(g_{i,j,k}\left(\textbf{X},t\right)\right)\right] \right)\right)\right) \textbf{f}_k \\
        \mathcal{L} = \sum_{k=1}^{m} \left( \sgn \left(\sum_{i=1}^{p_k} \left(\prod_{j=1}^{q_{p_k}} \left[1 - u\left(g_{i,j,k}\left(\textbf{X},t\right)\right)\right] \right)\right)\right) \mathcal{L}_k
    \end{gathered}
\end{equation}

The horizontally flipped unit step functions in Eq. \eqref{eqn:grashs_unitstep} are smoothed using sigmoid functions as follows (\cite{rashs}):

\begin{equation}
\label{eqn:rashs_sigmoid}
    \left[1 - u\left(g_{i,j,k}\left(\textbf{X},t\right)\right)\right] \approx \frac{1}{1 + e^{s \cdot g_{i,j,k}\left(\textbf{X},t\right)}}
\end{equation}

\begin{figure}[!h]
    \centering
    {\includegraphics[width=4in]{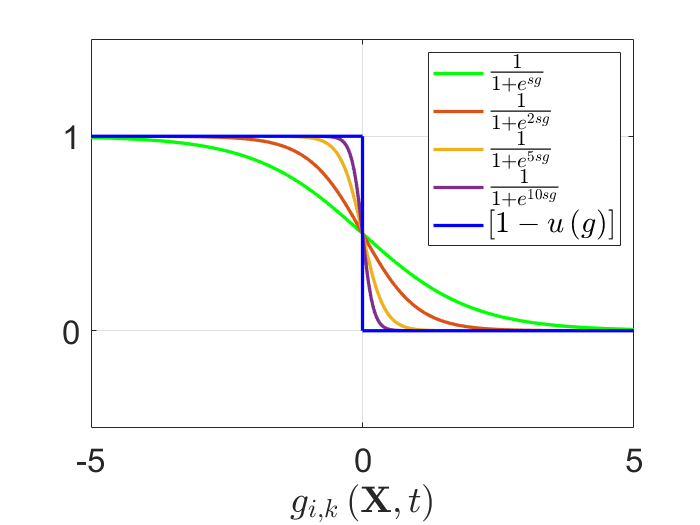}}
    \caption{\label{fig:unitStepApproximation}
    Approximation of unit step function using sigmoid function.}
\end{figure}

In Eq. \eqref{eqn:rashs_sigmoid}, $s$ is a measure of the slope at the transition point of the step. As $s \rightarrow \infty$, the sigmoid function approaches the horizontally flipped unit step function, as illustrated in Fig. \ref{fig:unitStepApproximation}.

\begin{figure}[!h]
    \centering
    {\includegraphics[width=4in]{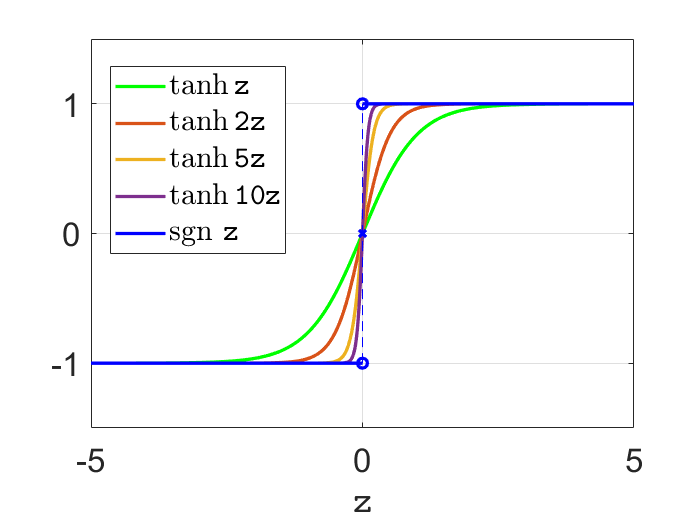}}
    \caption{\label{fig:signumApproximation}
    Approximation of signum function using hyperbolic tangent function.}
\end{figure}

The signum function in Eq. \eqref{eqn:grashs_unitstep}, which has a discontinuity at $0$, is smoothed using a hyperbolic tangent function:

\begin{equation}
    \label{eqn:grashs_tanh}
    \sgn\left(\texttt{z}\right) \approx \tanh{\zeta \texttt{z}}
\end{equation}

\noindent where $\zeta$ represents a parameter that represents the slope of the hyperbolic tangent function when $z=0$. As $\zeta \rightarrow \infty$, the hyperbolic tangent function approaches the signum function, as illustrated in Fig. \ref{fig:signumApproximation}.

Note that the saturation of the summation representing the OR logic can also be accomplished using a unit step function. However, the signum function is favored because its approximation using the hyperbolic tangent function also evaluates to $0$ when the input (summation of the minterms) is $0$, as observed in Fig. \ref{fig:signumApproximation}. This does not hold true for the sigmoid function approximating the unit step function, which evaluates to $0.5$ when the input is $0$, as observed in Fig. \ref{fig:unitStepApproximation}. This in turn introduces a significant amount of error in the smoothed equations of motion and Lagrangian. After applying the signum and hyperbolic tangent functions, Eq. \ref{eqn:grashs_unitstep} becomes:

\begin{equation}
    \label{eqn:grashs_final}
    \begin{gathered}
        \dot{\textbf{X}} = \sum_{k=1}^{m} \left( \tanh \left(\zeta \sum_{i=1}^{p_k} \left(\prod_{j=1}^{q_{p_k}} \left[\frac{1}{1 + e^{s \cdot g_{i,k}\left(\textbf{X},t\right)}}\right] \right)\right)\right) \textbf{f}_k \\
        \mathcal{L} = \sum_{k=1}^{m} \left( \tanh \left(\zeta \sum_{i=1}^{p_k} \left(\prod_{j=1}^{q_{p_k}} \left[\frac{1}{1 + e^{s \cdot g_{i,k}\left(\textbf{X},t\right)}}\right] \right)\right)\right) \mathcal{L}_k
    \end{gathered}
\end{equation}

Eq. \eqref{eqn:grashs_final} is a smooth approximation of the piecewise equations of motion and Lagrangian defined in Eq. \eqref{eqn:rashs_singleExp} when $B_k$ is represented using an arbitrary boolean logic. The boolean logic is essentially embedded into the equations of motion and Lagrangian and smoothed. This reduces the MPBVP that constitutes the necessary conditions of optimality in the indirect trajectory optimization framework to a simpler TPBVP because the associated interior-point boundary conditions are implicitly accounted for by Eq. \eqref{eqn:grashs_final}. Through homotopy, the quantities $s$ and $\zeta$ can be incrementally set to arbitrarily large values to make the TPBVP an arbitrarily close approximation of the MPBVP.

It must also be noted that if the DNF representation of $B_k$ consists of only one minterm for every flight segment, Eq. \eqref{eqn:grashs_final} collapses to:

\begin{equation}
    \label{eqn:rashs_final}
    \begin{gathered}
        \dot{\textbf{X}} = \sum_{k=1}^{m} \left(\prod_{j=1}^{q_k} \left[\frac{1}{1 + e^{s \cdot g_{j,k}\left(\textbf{X},t\right)}} \right]\right) \textbf{f}_k \\
        \mathcal{L} = \sum_{k=1}^{m} \left(\prod_{i=j}^{q_k} \left[\frac{1}{1 + e^{s \cdot g_{j,k}\left(\textbf{X},t\right)}} \right]\right) \mathcal{L}_k
    \end{gathered}
\end{equation}

\noindent which is essentially the RASHS formulation (\cite{rashs}) consisting of only AND logic.

\section{Trajectory Optimization using Indirect Methods}
\label{sec:trajectoryOptimization}
This investigation follows the indirect methods of trajectory optimization process, which involves solving the necessary conditions of optimality. For multi-phase systems, these conditions are represented by an MPBVP in a system of DAEs. When following the GRASHS approach, the MPBVP collapses to a TPBVP. The DAEs are solved using finite difference methods and homotopy. The necessary conditions of optimality and homotopy method are covered in Sections \ref{sec:necessaryConditions} and \ref{sec:continuation} respectively.

\subsection{Necessary Conditions of Optimality}
\label{sec:necessaryConditions}
Given a multi-phase trajectory consisting of $m$ flight segments, the equations of motion for segment $k$ is:

\begin{equation}
    \begin{gathered}
        \dot{\textbf{X}}=\textbf{f}_k\left(\textbf{X},\textbf{U},t\right) \text{, } t \in \left[t_{k-1}, ~ t_k\right] \\
        k=1, ~ 2, ~ ... ~ m
    \end{gathered}
\end{equation}

\noindent where $\textbf{U}$ is the control vector. The cost functional to be minimized is as follows:

\begin{equation}
    J=\textbf{$\Phi$}\left(\textbf{X}\left(t_m\right),t_m\right) + \sum_{k=1}^{m}\int\limits_{t_{k-1}}^{t_k} \mathcal{L}_k(\textbf{X},\textbf{U},t)dt
\end{equation}

\noindent where $\Phi$ is the terminal cost and the integral containing $\mathcal{L}$ is the path cost. The time $t_m$ corresponding to the end of the final segment $m$ is also the final time $t_f$ of the overall trajectory.

In addition to the equations of motion, the vehicle is also subject to the following end-point and interior-point boundary conditions:

\begin{equation}
    \label{eqn:allBoundaryConditions}
    \begin{gathered}
        \bm{\Psi}_0\left(\textbf{X}\left(t_0\right),t_0\right) = \textbf{0} \\
        \bm{\Psi}_k\left(\textbf{X}\left(t_k\right),t_k\right) = \textbf{0} \\
        \text{where } k = 1, ~ 2, ~ ..., ~ m
    \end{gathered}
\end{equation}

The interior-point boundary conditions are represented by $\bm{\Psi}_k = \textbf{0}$ for $k=1..,m-1$. These consist of the predicates that trigger the transition of flight segments and equality constraints to enforce continuity of the continuous states. The necessary conditions of optimality for this problem is an MPBVP in a system of DAEs. In this investigation, the solution of the GRASHS approach will be compared against that of the MPBVP. The derivation of the MPBVP when the switching conditions are represented as an arbitrary discrete logic is complicated and beyond the scope of this investigation. Instead, the MPBVP used for comparison will eliminate the OR logic by assuming it is known apriori which minterm from Eq. \eqref{eqn:grashs_Bk} causes $B_k$ to transition from $\texttt{false}$ to $\texttt{true}$. This apriori information will be obtained from the solution of the GRASHS approach. Assuming that the $i_k$\textsuperscript{th} minterm causes $B_k$ to transition to $\texttt{true}$, the interior-point boundary conditions are as follows:

\begin{equation}
    \begin{gathered}
        \bm{\Psi}_k\left(\textbf{X}_k,t_k\right) = 
        \begin{pmatrix}
          \bm{X}\left(t_k^-\right) - \bm{X}\left(t_k^+\right) \\
          \prod\limits_{j=1}^{q_{p_k}} g_{i_k,j,k}\left(\textbf{X},t\right)
        \end{pmatrix}
        =
        \textbf{0} \\
        \text{where } k = 1, ~ 2, ~ ..., ~ m-1
    \end{gathered}
\end{equation}

Accordingly, the MPBVP representing the necessary conditions of optimality is given as follows (\cite{brysonho}):

\begin{equation}
    \label{eqn:mpbvp_necessary_conditions}
    \begin{gathered}
        \dot{\bm{\lambda}}=-\left(\frac{\partial \textbf{H}_k}{\partial \textbf{X}}\right)^T \text{, } t \in \left[t_{k-1}, ~ t_k\right] \\
        \dot{\textbf{X}}=\textbf{f}_k\left(\textbf{X},\textbf{U},t\right)\text{, } t \in \left[t_{k-1}, ~ t_k\right] \\
        \frac{\partial \textbf{H}_k}{\partial \textbf{U}} = \textbf{0} \left[t_{k-1}, ~ t_k\right] \\
        \bm{\Psi}_0\left(\textbf{X}\left(t_0\right),t_0\right) = \textbf{0} \\
        \bm{\Psi}_k\left(\textbf{X}\left(t_k\right),t_k\right) = \textbf{0} \\
        \textbf{H}_{m}\left(t_m\right) dt_m - \textbf{H}_{1}\left(t_0\right) dt_0 - \bm{\lambda}^T\left(t_m\right) d\textbf{X}\left(t_m\right) + \bm{\lambda}^T\left(t_0\right) d\textbf{X}\left(t_0\right) + \\
        d{\Phi}\left(\textbf{X}\left(t_m\right),t_m\right) = 0 \\
        {\bm{\lambda}^T}\left(t^-_{k\neq m}\right) = {\bm{\lambda}^T}\left(t^+_{k\neq m}\right) + \bm{\Pi}_k^T \frac{\partial \bm{\Psi}_k}{\partial \textbf{X}\left(t_k\right)} \\
        \textbf{H}_{k-1}\left(t_k^-\right) = \textbf{H}_k\left(t_k^+\right) - \bm{\Pi}_k^T \frac{\partial \bm{\Psi}_k}{\partial t_k} \\
        \text{where } k=1, ~ 2, ~...,~ m  \text{ and } \textbf{H}_k = \mathcal{L}_k + \bm{\lambda}^T \textbf{f}_k
    \end{gathered}
\end{equation}

\noindent where $\textbf{H}$ is the Hamiltonian and $\bm{\lambda}$ is the co-state vector. This MPBVP must be solved numerically using methods such as multiple shooting (\cite{shooting}) and finite difference (\cite{collocation}). This MPBVP is difficult to solve because the numerical methods require an initial guess for each flight segment, which is not straightforward. The challenge is exacerbated as the number of flight segments increases because the number of interior-point boundary conditions that must be enforced also increases. The GRASHS approach addresses the latter issue because these interior-point boundary conditions are embedded into the equations of motion in Eq. \ref{eqn:grashs_final} and are not required to be explicitly enforced. Therefore, the necessary conditions of optimality collapse to a TPBVP as follows (\cite{brysonhoTPBVP}):

\begin{equation}
    \label{eqn:tpbvp_necessary_conditions}
    \begin{gathered}
        \dot{\textbf{X}}=\textbf{f}\left(\textbf{X},\textbf{U},t\right) \\
        \dot{\textbf{$\lambda$}}=-\left(\frac{\partial \textbf{H}}{\partial \textbf{X}}\right)^T \\
        \frac{\partial \textbf{H}}{\partial \textbf{U}} = \textbf{0} \\
        \bm{\Psi}_0\left(\textbf{X}\left(t_0\right),t_0\right) = \textbf{0} \\
        \bm{\Psi}_f\left(\textbf{X}\left(t_f\right),t_f\right) = \textbf{0} \\
        \textbf{H}\left(t_f\right) dt_f - \textbf{H}\left(t_0\right) dt_0 - \bm{\lambda}^T\left(t_f\right) d\textbf{X}\left(t_f\right) + \bm{\lambda}^T\left(t_0\right) d\textbf{X}\left(t_0\right) + \\
        d\Phi\left(\textbf{X}\left(t_f\right),t_f\right) = 0 
    \end{gathered}
\end{equation}

\noindent where the subscript $f$ represents the quantity at final time $t_f$. The generation of an initial guess to solve Eq. \eqref{eqn:tpbvp_necessary_conditions} to guarantee convergence of the numerical methods to a solution is still a challenging task. This challenge is mitigated using homotopy, as explained in Section \ref{sec:continuation}.

\subsection{Numerical Solution to Necessary Conditions using Homotopy}
\label{sec:continuation}
In the homotopy process (\cite{continuation}), rather than directly attempting to solve the original TPBVP in Eq. \eqref{eqn:tpbvp_necessary_conditions}, a trivially simple problem with a short time of flight is solved. This problem must be trivial enough to enable convergence to a solution even with a poor initial guess. After the trivial problem is solved, the problem is evolved in steps to the original problem of interest by gradually modifying the boundary conditions. If the problem in a given step varies only slightly from the preceding step, the solution from the preceding step will be close to that of the current step. Therefore, it can be used as the guess for the current step, thereby improving convergence when compared to attempting to solve the original TPBVP outright. Homotopy has been successfully applied to solve a variety of trajectory optimization problems using indirect methods (\cite{rashs, continuation2, continuation3, continuation4}).

When employing the GRASHS approach, in addition to the boundary conditions, the homotopy process is also embedded with $s$ and $\zeta$ that control the slope of the sigmoid and hyperbolic tangent functions at the flight segment transition points. These quantities are initially seeded with low values for a gradual transition in $\textbf{f}$ and $\mathcal{L}$, thereby further improving convergence. These parameters are then increased to arbitrarily large values through homotopy to bring the solution arbitrarily close to that of the MPBVP.

The homotopy process does involve trial and error with reference to the choice of homotopy parameters, the number of iterations and the degree of variation of the parameters between iterations. Recent advances in the indirect framework have shown to mitigate this through the use of adaptive continuation (\cite{adaptiveContinuation}) and multistage stabilized continuation (\cite{multiStageContinuation}). However, these advancements are not employed here and are not the focus of this investigation.

\section{Flight Dynamics Model}
\label{sec:eom}
This investigation assumes the 3 degrees-of-freedom (DOF) flight dynamics model as described in \cite{rashs}. This model assumes a spherical rotating planet of uniform mass density whose center is assumed to be inertial. Accordingly, the state variables are altitude ($h$), longitude ($\theta$), latitude ($\phi$), atmospheric-relative velocity ($v$), atmospheric-relative flight path angle ($\gamma$), heading angle ($\psi$) and the mass of consumed fuel ($m_F$). The equations of motion are as follows:

\begin{equation}
    \label{eqn:eom_edl}
    \begin{gathered}
        \frac{^id \textbf{r}}{d t} = \left(\bm{\omega} \times \textbf{r}\right) + \textbf{v} \\
        \frac{^id}{dt}\left(\left(\bm{\omega} \times \textbf{r}\right) + \textbf{v}\right) = \frac{\textbf{F}}{m} \\
        \frac{d m_F}{d t} = \dot{m}_{F, max} \frac{T}{T_{max}}
    \end{gathered}
\end{equation}

\noindent where:
	
\begin{equation}
    \label{eqn:inertialPositionAndVelocity}
    \begin{gathered}
        \textbf{r} = \left(R + h\right) \left( \cos{\phi} \cos{\theta} ~ \bm{\hat{\textbf{X}}}_G + \cos{\phi} \sin{\theta} ~ \bm{\hat{\textbf{Y}}}_G + \sin{\phi} ~\bm{\hat{\textbf{Z}}}_G\right) \\
        \textbf{v} = v \left(\cos{\gamma} \sin{\psi} ~ \bm{\hat{\textbf{e}}}_E + \cos{\gamma} \cos{\psi} ~ \bm{\hat{\textbf{e}}}_N + \sin{\gamma} ~ \bm{\hat{\textbf{e}}}_Z\right)
    \end{gathered}
\end{equation}

\begin{figure}[!t]
	\centering
	\begin{subfigure}[t]{0.45\textwidth}
		\centering
		\includegraphics[width=\textwidth]{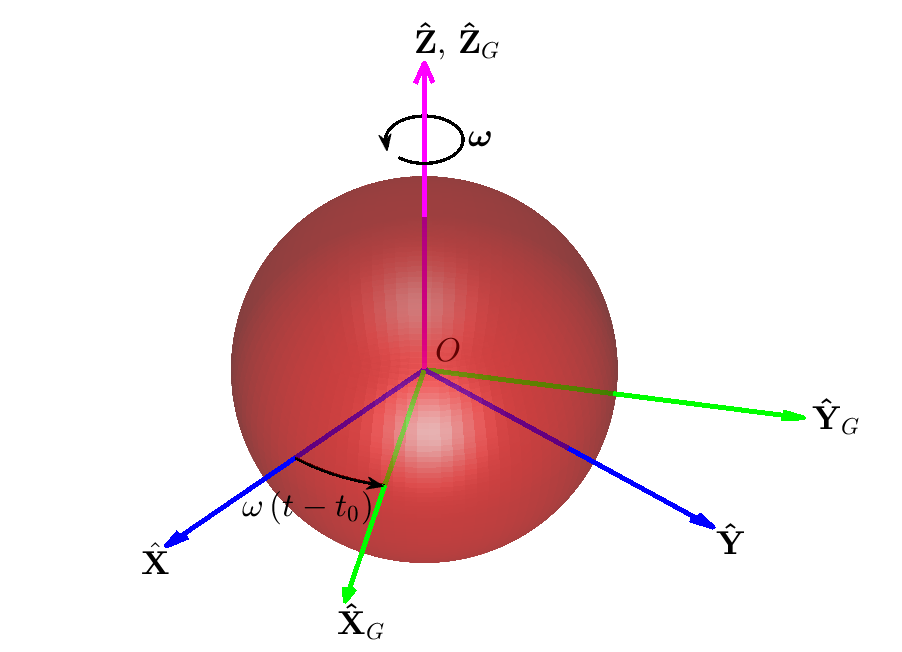}
		\caption{\label{fig:inertialToGeographic}Relationship between inertial frame and PCPF. Reprinted from \cite{rashs}.}
	\end{subfigure}
	\begin{subfigure}[t]{0.45\textwidth}
		\centering
		\includegraphics[width=\textwidth]{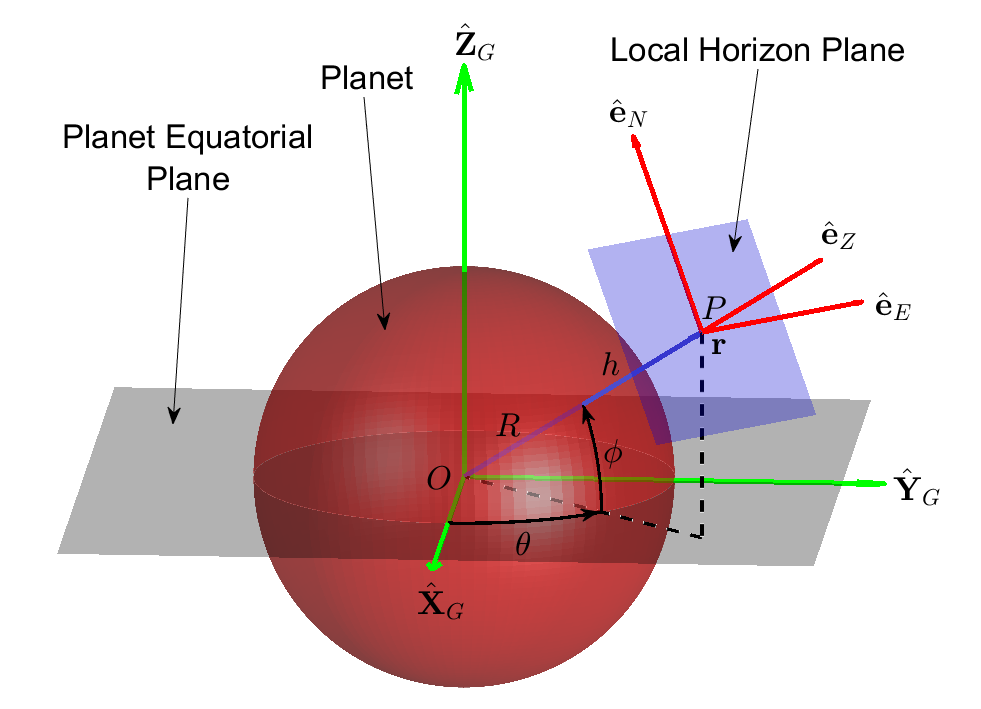}
		\caption{\label{fig:geographicToHorizon}Relationship between PCPF and local horizon frame. Reprinted from \cite{rashs}.}
	\end{subfigure}
	\begin{subfigure}[t]{0.45\textwidth}
		\centering
		\includegraphics[width=\textwidth]{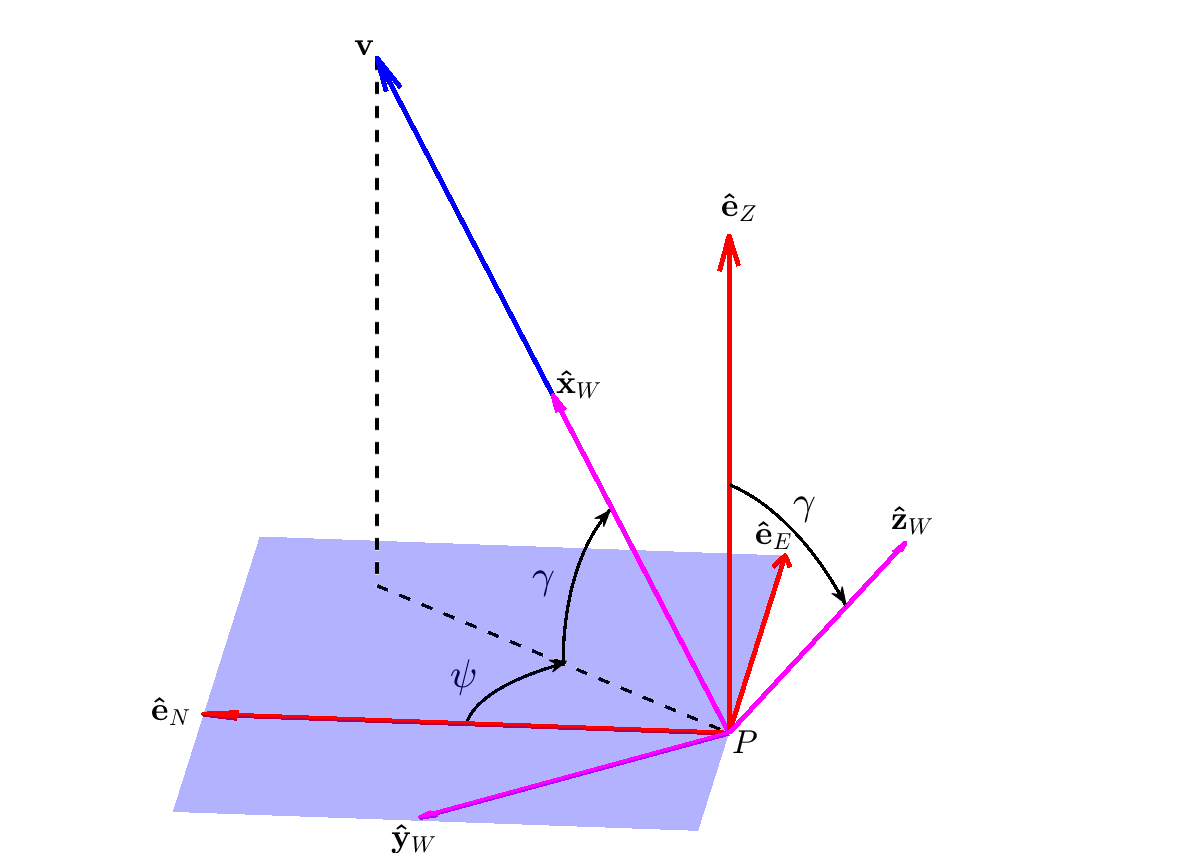}
		\caption{\label{fig:localHorizon}Relationship between local horizon and wind frames. Reprinted from \cite{rashs}.}
	\end{subfigure}
	\begin{subfigure}[t]{0.45\textwidth}
		\centering
		\includegraphics[width=\textwidth]{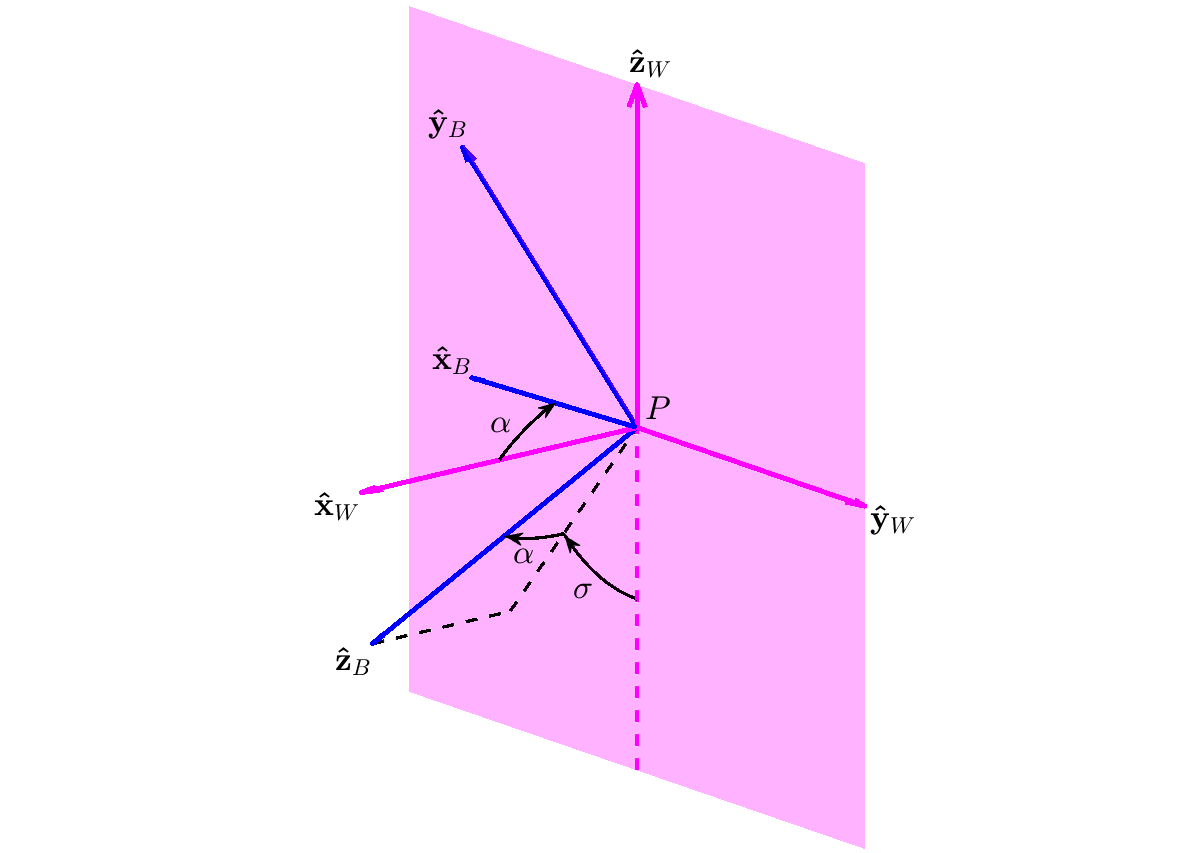}
		\caption{\label{fig:windFrame}Relationship between wind and body frames. Reprinted from \cite{rashs}.}
	\end{subfigure}
	\caption{\label{fig:coordinateFrames}Relationship between coordinate frames. Adopted from \cite{rashs}.}
\end{figure}

In Eq. \eqref{eqn:eom_edl}, $\textbf{F}$ is the force vector acting on the vehicle given by:

\begin{equation}
    \label{eqn:edl_force}
    \textbf{F} = \left(T \cos\alpha - D\right) ~ \bm{\hat{\textbf{x}}}_W + \left(T\sin\alpha + L\right) \sin{\sigma} ~ \bm{\hat{\textbf{y}}}_W +
    \left(T\sin\alpha + L\right) \cos{\sigma} ~ \bm{\hat{\textbf{z}}}_W - \frac{\mu m}{\left(R+h\right)^2} ~ \bm{\hat{\textbf{e}}}_Z
\end{equation}

\noindent where $T$ is the thrust, $\alpha$ is the angle-of-attack, $\sigma$ is the bank angle, $\mu$ is the standard gravitational acceleration of the planet and $R$ is the planetary radius. Also, $\bm{\hat{\textbf{x}}}_W$, $\bm{\hat{\textbf{y}}}_W$ and $\bm{\hat{\textbf{z}}}_W$ are the unit vectors defining the wind frame, as illustrated in Fig. \ref{fig:windFrame}. In this figure, $\bm{\hat{\textbf{x}}}_B$, $\bm{\hat{\textbf{y}}}_B$ and $\bm{\hat{\textbf{z}}}_B$ are the unit vectors defining the body frame. The quantities $L$ and $D$ are the lift and drag forces given by:

\begin{equation}
    \begin{gathered}
        L = \frac{1}{2} \rho_o e^{-\left(\frac{h}{H}\right)} v^2 C_L S \\
        D = \frac{1}{2} \rho_o e^{-\left(\frac{h}{H}\right)} v^2 C_D S
    \end{gathered}
\end{equation}

\noindent where $\rho_o$ is the surface atmospheric density, $H$ is the atmospheric scale height, $C_L$ is the lift coefficient, $C_D$ is the drag coefficient and $S$ is the reference area. Additionally, $\dot{m}_{F,max}$ is the maximum fuel mass flow rate corresponding to the maximum thrust $T_{max}$. The instantaneous mass $m$ of the vehicle is given by:

\begin{equation}
    m = m_0 - m_F
\end{equation}

\noindent where $m_0$ is the total initial mass of the vehicle at the beginning of a given flight segment.

In Eq. \eqref{eqn:inertialPositionAndVelocity}, $\textbf{r}$ is the inertial position vector and $\textbf{v}$ is the atmospheric-relative velocity vector. Also, $\bm{\hat{\textbf{X}}}_G$, $\bm{\hat{\textbf{Y}}}_G$ and $\bm{\hat{\textbf{Z}}}_G$ are the unit vectors defining the planet-centered planet-fixed (PCPF) frame and $\bm{\hat{\textbf{e}}}_E$, $\bm{\hat{\textbf{e}}}_N$ and $\bm{\hat{\textbf{e}}}_Z$ are the unit vectors defining the local East-North-Up (ENU) frame. These quantities are illustrated in Fig. \ref{fig:inertialToGeographic}, \ref{fig:geographicToHorizon} and \ref{fig:localHorizon}. The equations of motion in Eq. \eqref{eqn:eom_edl} represent $\textbf{f}$ in Eqs. \eqref{eqn:mpbvp_necessary_conditions} and \eqref{eqn:tpbvp_necessary_conditions}.

\section{Numerical Examples: Application of GRASHS to Optimize Mars EDL Trajectory}
This section applies the GRASHS approach to optimize a Mars EDL trajectory. The advantage of the GRASHS workflow is illustrated by the parachute deployment event, which is configured to be triggered by either a velocity or an altitude condition. Consequently, the conditions applicable for parachute descent are represented as a combination of AND and OR logic. Optimal trajectories are generated for two mission profiles that are identical in every aspect except the parachute deployment altitude, which is set to a lower value in the first profile. Therefore, the parachute deployment event will be shown to be velocity-triggered for the first profile, and altitude-triggered for the second profile. The GRASH result is compared with the RASHS and MPBVP solutions for both profiles. It is important to note that the GRASH formulation bears no apriori knowledge of which trigger event (altitude or velocity) will be hit first, which is in fact a key advantage of this formulation when solving the aforementioned EDL trajectory optimization problem. Upon gaining the knowledge of the trigger event from the GRASHS solution, the trajectories are re-calculated using RASHS and MPBVP for comparison.

\subsection{Mission Architecture}
The mission architecture considered in this investigation is similar to that of the Mars Science Laboratory (MSL) mission (\cite{msl1,msl2,msl3,msl4}), with the following simplifying assumptions to focus on the GRASHS approach rather than the architecture:

\begin{enumerate}
  \item The heat shield is jettisoned simultaneously with parachute deployment.
  \item The powered descent continues all the way to touch down and does not employ a sky crane.
\end{enumerate}

\begin{figure}[!h]
    \centering
    {\includegraphics[width=4in]{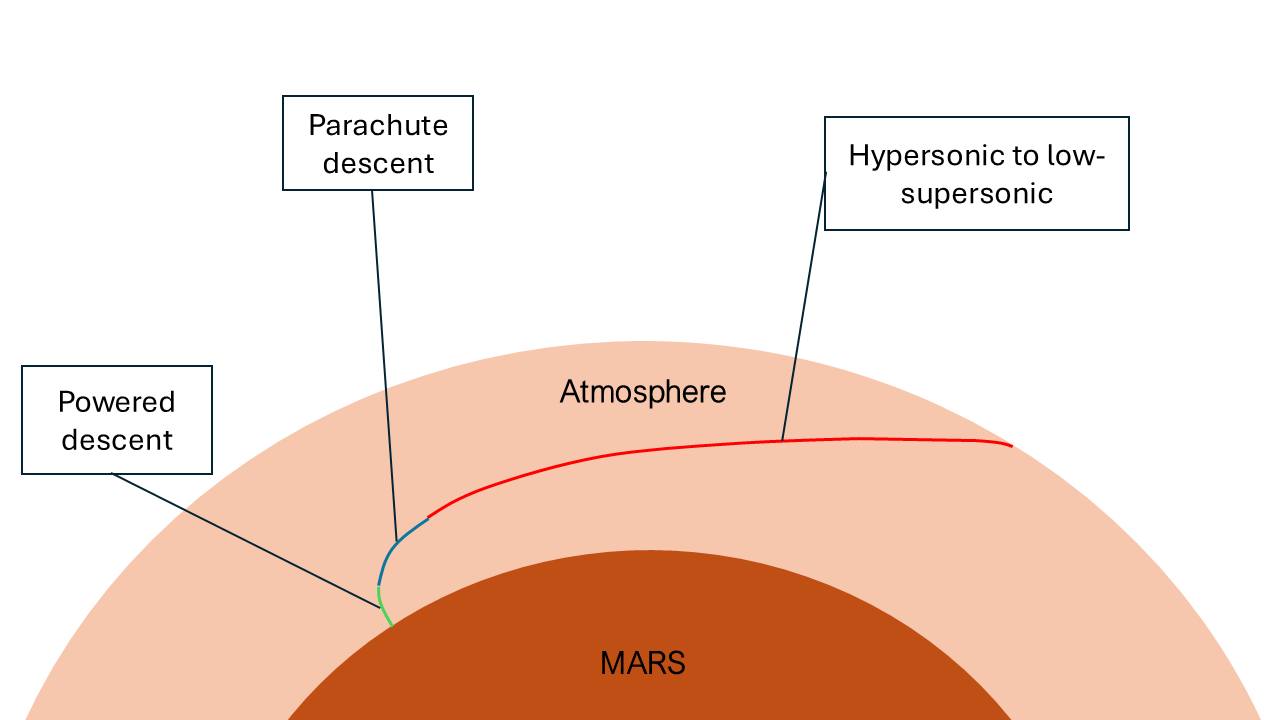}}
    \caption{\label{fig:mslMissionArch}
    Illustration of the mission architecture.}
\end{figure}

With these simplifications, the mission consists of the following flight segments as illustrated in Fig. \ref{fig:mslMissionArch}:

\begin{enumerate}
  \item \textbf{Segment 1 - hypersonic to low supersonic:} This segment begins at entry interface and ends at parachute deployment.
  \item \textbf{Segment 2 - parachute descent:} This segment begins at parachute deployment and ends at powered descent initiation.
  \item \textbf{Segment 3 - powered descent:} This segment begins at powered descent initiation and ends at touchdown.
\end{enumerate}

The Martian atmospheric entry interface is set at $h=120$ km, $v=5.9$ km/s, $\theta=0$ deg, $\phi=0$ deg and $\psi = 90$ deg. This marks the beginning of segment 1 (hypersonic to low supersonic). The total entry mass is $3152$ kg. The vehicle is trimmed at a nonzero $\alpha$ by means of a center-of-mass offset achieved using ballasts bearing a mass of 150 kg. This provides a constant $C_L$ of $0.25$ and a $C_D$ of $1.24$ for this segment. The reference area ($S$) of the vehicle is $15.9$ m\textsuperscript{2}. The vehicle is maneuvered by modulating the bank angle. This segment ends when the vehicle decelerates to the parachute deployment velocity ($v_P$) of $408$ m/s or descents to the parachute deployment altitude ($h_{P}$). The value of $h_P$ is set to $3.5$ km for mission profile 1 and $6.5$ km for mission profile 2. The end of segment 1 marks the beginning of segment 2 (parachute descent).

When the parachute descent segment begins, the vehicle jettisons the ballasts and the heat shield. The mass of the heat shield is assumed to be $385$ kg. Therefore, at parachute deployment, the vehicle mass reduces to $2617$ kg. Upon jettisoning the ballasts, the vehicle trims at $0$ deg $\alpha$, resulting in a $C_L$ of $0$. With the parachute deployed, $C_D$ changes to $9.43$. During parachute descent, the vehicle cannot be maneuvered. This segment ends when the vehicle descends below the powered descent initiation altitude ($h_{PDI}$) of $2$ km and segment 3 (powered descent) begins.

When the powered descent segment begins, the backshell and the parachute are jettisoned. These bring the vehicle mass further down to $2268$ kg and $C_D$ to $0.31$. Throughout this segment, $\alpha$ is held at $0$ deg with no sideslip. Therefore, $C_L$ remains at $0$ and the thrust vector is always oriented retrograde. The vehicle mass trends down throughout this segment, consistent with the propellant mass flow rate ($\dot{m}_F$), which in turn governs the amount of thrust produced. The descent engines bear a collective $I_{sp}$ of $210$ s and a maximum propellant mass flow rate of $12.43$ kg/s, which in turn translates to a maximum thrust of $25.6$ kN. The total propellant onboard at powered descent initiation is $387$ kg. This segment ends at touchdown. At touchdown, $h=0$ km, $v = 0.1$ m/s, $\theta = 16.027$ deg and $\phi = 1.1809$ deg. 

Table \ref{tab:flightSegmentsControl} summarizes the control mechanism. Table \ref{tab:flightSegmentsMassAero} summarizes the mass, aerodynamic characteristics and maximum propellant mass flow rate for each flight segment. Table \ref{tab:flightSegmentConditions} summarizes the conditions that determine the active flight segment. Note that the conditions for Segment 2 are composed of a combination of AND and OR logic. Table \ref{tab:boundaryConditions} summarizes the initial and final conditions on the state variables.

\begin{table}[!h]
    \centering
    \caption{\label{tab:flightSegmentsControl}Control mechanism for each flight segment.}
    \begin{tabular}{lc}
        \hline
        Flight Segment & Control Mechanism \\\hline
        Segment 1 & Bank angle modulation \\
        Segment 2 & None \\
        Segment 3 & Retrograde thrust \\
        \hline
    \end{tabular}
\end{table}

\begin{table}[!h]
    \centering
    \caption{\label{tab:flightSegmentsMassAero}Mass, aerodynamic characteristics and maximum propellant mass flow rate for each flight segment.}
    \begin{tabular}{lccccc}
        \hline
        Flight Segment & Mass $m$, kg & $C_L$ & $C_D$ & $\dot{m}_{F,max}$, kg/s \\\hline
        Segment 1 & $3152$ & $0.25$ & $1.24$ & 0 \\
        Segment 2 & $2617$ & $0$ & $9.43$ & 0 \\
        Segment 3 & $2268$ (initial) & $0$ & $0.31$ & 12.43 \\
        \hline
    \end{tabular}
\end{table}

\begin{table}[!h]
    \centering
    \caption{\label{tab:flightSegmentConditions}Conditions that must be satisfied for each flight segment.}
    \begin{tabular}{lc}
        \hline
        Flight Segment & Conditions \\\hline
        Segment 1 & $\left(v \geq v_P\right)  \land \left(h \geq h_P\right)$ \\
        Segment 2 & $\left(\left(v < v_P\right)  \lor \left(h < h_P\right)\right) \land \left(h \geq h_{PDI}\right)$ \\
        Segment 3 & $h < h_{PDI}$ \\
        \hline
    \end{tabular}
\end{table}

\begin{table}[!h]
    \centering
    \caption{\label{tab:boundaryConditions}Initial and final conditions.}
    \begin{tabular}{lcc}
        \hline
        State & Initial condition & Final condition \\\hline
        Altitude $h$ & $120$ km & $0$ km \\
        Longitude $\theta$ & $0$ deg & $16.027$ deg \\
        Latitude $\phi$ & $0$ deg & $1.1809$ deg \\
        Atmospheric-relative velocity $v$ & $5.9$ km/s & $0.1$ m/s \\
        Atmospheric-relative flight-path angle $\gamma$ & Free & Free \\
        Atmospheric-relative heading angle $\psi$ & $90$ deg & Free \\
        Mass of propellant consumed $m_F$ & $0$ kg & Free \\
        \hline
    \end{tabular}
\end{table}

The surface density ($\rho_o$) and scale height ($H$) of the Martian atmosphere are assumed to be $0.025$ kg/m\textsuperscript{3} and $11.1$ km, respectively.

The stagnation-point heat-load on the vehicle must be minimized from entry interface to powered descent initiation. The stagnation-point heat-load, $Q$, is calculated by integrating the stagnation-point heat-rate, $\dot{q}$, given by \cite{suttonGraves}, over time:

\begin{equation}
    Q\left(t\right) = \int_0^{t} \dot{q} d \tau = \int_0^{t} k \sqrt{\frac{\rho_0 e^{-\frac{h}{H}}}{R_N}} v^3 d\tau
\end{equation}

\noindent where $k$ is an empirical constant and $R_N$ is the nose radius of the vehicle. For the vehicle under consideration in this mission, $R_N = 1.125$ m.

During powered descent, the thrust must be minimized, which acts as a surrogate for minimizing fuel consumption. Accordingly, the cost functional for the overall mission is:

\begin{equation}
    \label{eqn:missionCost}
    J = K_1 \int_0^{t_P} \sqrt{\frac{\rho_o e^{-\frac{h}{H}}}{R_N}} v^3 d\tau +
    K_2 \int_{t_P}^{t_{PDI}} \sqrt{\frac{\rho_o e^{-\frac{h}{H}}}{R_N}} v^3 d\tau +
    K_3 \int_{t_{PDI}}^{t_f} T^2 d\tau
\end{equation}

\noindent where $t_P$ and $t_{PDI}$ are the times at parachute deployment and PDI respectively, and $K_1$, $K_2$ and $K_3$ are weights. The resultant $\mathcal{L}$ for each flight segment is summarized in Table \ref{tab:missionLagrangian}.

\begin{table}[!h]
    \centering
    \caption{\label{tab:missionLagrangian}Lagrangian for each flight segment.}
    \begin{tabular}{lc}
        \hline
        Flight Segment & Lagrangian $\mathcal{L}$ \\\hline
        Segment 1 & $\mathcal{L}_1 = K_1 \sqrt{\frac{\rho_o e^{-\frac{h}{H}}}{R_N}} v^3$ \\
        Segment 2 & $\mathcal{L}_2 = K_2 \sqrt{\frac{\rho_o e^{-\frac{h}{H}}}{R_N}} v^3$ \\
        Segment 3 & $\mathcal{L}_3 = K_3 T^2$ \\
        \hline
    \end{tabular}
\end{table}

In this problem, the weights $K_1$ and $K_2$ are set to $1$ and $K_3$ is set to $10^{-5}$ to bring the sum of path costs for Segment 1 and Segment 2 to roughly the same order of magnitude as the path cost for Segment 3. It is required to calculate a trajectory with flight segment transitions governed by the conditions defined in Table \ref{tab:flightSegmentConditions}, such that it minimizes the cost functional in Eq. \eqref{eqn:missionCost}, subject to the end-point boundary conditions defined in Table \ref{tab:boundaryConditions}. As mentioned, two mission profiles are considered, where $h_P = 3.5$ km  for profile 1 and $h_P = 6.5$ km for profile 2.

\subsection{Mission Profile 1: Low Parachute Deployment Altitude}
\label{sec:missionProfile1}
This section demonstrates the optimal trajectory generation using the GRASHS approach when the parachute deployment altitude, $h_P$, is $3.5$ km. To aid in the solution process, the states $h$ and $v$ are scaled by their initial conditions ($h\left(0\right) = 120$ km and $v\left(0\right) = 5.9$ km/s), and the state $m_F$ is scaled by the total propellant onboard prior to PDI ($387$ kg). This results in the corresponding scaled variables $\hbar$, $V$ and $M_F$, such that:

\begin{equation}
    \begin{gathered}
        \hbar = \frac{h}{h\left(0\right)} \\
        V = \frac{v}{v\left(0\right)} \\
        M_F = \frac{m_F}{387}
    \end{gathered}
\end{equation}

Therefore, the state vector of the scaled problem is $\textbf{X} = \left[\hbar ~ \theta ~ \phi ~ V ~ \gamma ~ \psi ~ M_F\right]^T$. The DNF representations of the conditions defined in Table \ref{tab:flightSegmentConditions} are:

\begin{equation}
    \label{eqn:DNF_GRASHS}
    \begin{gathered}
        \textrm{Segment 1: } \left(\left(\frac{v_P}{v\left(0\right)} - V\right) < 0\right) \cdot \left(\left(\frac{h_P}{h\left(0\right)} - \hbar\right) < 0 \right) \\
        \textrm{Segment 2: } \left(\left(V - \frac{v_P}{v\left(0\right)}\right) < 0\right) \cdot \left(\left(\frac{h_{PDI}}{h\left(0\right)} - \hbar\right) < 0 \right) + \\
        \left(\left(\hbar - \frac{h_P}{h\left(0\right)}\right) < 0\right) \cdot \left(\left(\frac{h_{PDI}}{h\left(0\right)} - \hbar\right) < 0 \right) \\
        \textrm{Segment 3: } \left(\left(\hbar - \frac{h_{PDI}}{h\left(0\right)}\right) < 0\right)
    \end{gathered}
\end{equation}

\noindent with the following predicates:

\begin{equation}
    \label{eqn:predicates}
    \begin{gathered}
        g_{1,1,1} < 0 \Longrightarrow \left(\frac{v_P}{v\left(0\right)} - V\right) < 0 \\
        g_{1,2,1} < 0 \Longrightarrow \left(\frac{h_P}{h\left(0\right)} - \hbar\right) < 0 \\
        g_{1,1,2} < 0 \Longrightarrow \left(V - \frac{v_P}{v\left(0\right)}\right) < 0 \\
        g_{1,2,2} < 0 \Longrightarrow \left(\frac{h_{PDI}}{h\left(0\right)} - \hbar\right) < 0 \\
        g_{2,1,2} < 0 \Longrightarrow \left(\hbar - \frac{h_P}{h\left(0\right)}\right) < 0 \\
        g_{2,2,2} < 0 \Longrightarrow \left(\frac{h_{PDI}}{h\left(0\right)} - \hbar\right) < 0 \\
        g_{1,1,3} < 0 \Longrightarrow \left(\hbar - \frac{h_{PDI}}{h\left(0\right)}\right) < 0
    \end{gathered}
\end{equation}

It is worth mentioning that in Eq. \eqref{eqn:predicates}, $\left(g_{1,1,2} < 0\right)$ is the NOT of $\left(g_{1,1,1} < 0\right)$, $\left(g_{2,1,2} < 0\right)$ if the NOT of $\left(g_{1,2,1} < 0\right)$, and $\left(g_{1,1,3} < 0\right)$ is the NOT of $\left(g_{1,2,2} < 0\right)$ and $\left(g_{2,2,2} < 0\right)$. By applying the GRASHS approach, the continuous and differentiable equations of motion that embed the predicates summarized in Eq. \eqref{eqn:predicates} are given by:

\begin{multline}
    \label{eqn:eom_GRASHS}
    \textbf{f} = \left[\left(\frac{1}{1 + e^{s \left(\frac{v_P}{v\left(0\right)} - V\right)}}\right)
    \left(\frac{1}{1 + e^{s \left(\frac{h_P}{h\left(0\right)} - \hbar\right)}}\right)\right] \textbf{f}_1 + \\
    \tanh\left[\zeta\left(\left(\left(\frac{1}{1 + e^{s \left(V - \frac{v_P}{v\left(0\right)}\right)}}\right)
    \left(\frac{1}{1 + e^{s \left(\frac{h_{PDI}}{h\left(0\right)} - \hbar\right)}}\right)\right) +
    \left(\left(\frac{1}{1 + e^{s \left(\hbar - \frac{h_P}{h\left(0\right)}\right)}}\right)
    \left(\frac{1}{1 + e^{s \left(\frac{h_{PDI}}{h\left(0\right)} - \hbar\right)}}\right)\right)\right)\right] \textbf{f}_2 + \\
    \left[\left(\frac{1}{1 + e^{s \left(\hbar - \frac{h_{PDI}}{h\left(0\right)} \right)}}\right)\right] \textbf{f}_3
\end{multline}

\noindent where $\textbf{f}_1$, $\textbf{f}_2$ and $\textbf{f}_3$ are the equations of motion for the flight segments 1, 2 and 3, respectively. Note that the equations are the same as those defined in Section \ref{sec:eom}, but with different mass, aerodynamic coefficients and maximum propellant mass flow rate for each segment, as defined in Table \ref{tab:flightSegmentsMassAero}. This makes the original equations of motion piecewise continuous, which is made continuous and differential by GRASHS in Eq. \eqref{eqn:eom_GRASHS}.

Similarly, applying GRASHS on the Lagrangian yields:
\begin{multline}
    \label{eqn:lagrangian_GRASHS}
    \mathcal{L} = \left[\left(\frac{1}{1 + e^{s \left(\frac{v_P}{v\left(0\right)} - V\right)}}\right)
    \left(\frac{1}{1 + e^{s \left(\frac{h_P}{h\left(0\right)} - \hbar\right)}}\right)\right] \mathcal{L}_1 + \\
    \tanh\left[\zeta\left(\left(\left(\frac{1}{1 + e^{s \left(V - \frac{v_P}{v\left(0\right)}\right)}}\right)
    \left(\frac{1}{1 + e^{s \left(\frac{h_{PDI}}{h\left(0\right)} - \hbar\right)}}\right)\right) +
    \left(\left(\frac{1}{1 + e^{s \left(\hbar - \frac{h_P}{h\left(0\right)}\right)}}\right)
    \left(\frac{1}{1 + e^{s \left(\frac{h_{PDI}}{h\left(0\right)} - \hbar\right)}}\right)\right)\right)\right] \mathcal{L}_2 + \\
    \left[\left(\frac{1}{1 + e^{s \left(\hbar - \frac{h_{PDI}}{h\left(0\right)} \right)}}\right)\right] \mathcal{L}_3
\end{multline}

\noindent where $\mathcal{L}_1$, $\mathcal{L}_2$ and $\mathcal{L}_3$ are the Lagrangians for each flight segment, as defined in Table \ref{tab:missionLagrangian}. Eqs. \eqref{eqn:eom_GRASHS} and \eqref{eqn:lagrangian_GRASHS} and the boundary conditions defined in Table \ref{tab:boundaryConditions} are plugged into the TPBVP in Eq. \eqref{eqn:tpbvp_necessary_conditions}, whose solution yields the optimal multi-phase EDL trajectory.

The TPBVP is solved using homotopy, which is implemented in five steps, each consisting $500$ or $1000$ iterations to evolve a trivial optimization problem towards the desired problem. The TPBVP in each iteration is solved using \textit{bvp4c} (\cite{bvp4c1, bvp4c2}), which is a MATLAB implementation of finite difference methods to solve a system of DAEs. To seed the homotopy process, the DAEs in the TPBVP are propagated for $0.1$ seconds using \textit{ode45} (\cite{ode45}), which is a MATLAB implementation of the Dormand-Prince method (\cite{dormandPrince}). The initial conditions for the propagation are as follows:

\begin{equation}
    \begin{gathered}
        \textbf{X}^T = 
        \begin{bmatrix}
            \hbar & \theta & \phi & V & \gamma & \psi & m_F
        \end{bmatrix}
        =
        \begin{bmatrix}
            1 & 0 & 0 & 1 & 0 & \frac{\pi}{2} & 0
        \end{bmatrix}
        \\
        \bm{\lambda}^T =
        \begin{bmatrix}
            \lambda_{\hbar} & \lambda_{\theta} & \lambda_{\phi} & \lambda_V & \lambda_{\gamma} & \lambda_{\psi} & \lambda_{m_F}
        \end{bmatrix}
        =
        \begin{bmatrix}
            0 & 1 & -1 & -1 & 0 & 0 & 0
        \end{bmatrix}
    \end{gathered}
\end{equation}

Additionally, for this propagation, $C_L$ for Segment 1 is set to $0$, and the slope parameters $s$ and $\zeta$ in the sigmoid and the hyperbolic tangent functions are set to $100$ and $1$, respectively. Through homotopy, $C_L$ for Segment 1 will eventually be brought to the desired value of $0.25$, and $s$ and $\zeta$ will be brought up to $40,000$.

The five homotopy steps are implemented as follows:
\begin{enumerate}
    \item The first homotopy step iterates over the final boundary conditions on $\theta$ and $\phi$. The final boundary condition on $v$ and the initial boundary condition on $\psi$ are set free. For the first iteration, the final boundary conditions on $\theta$ and $\phi$ are set to the corresponding final values from the \textit{ode45} propagation. This iteration constitutes the trivial problem, and the \textit{ode45} propagation result serves as the initial guess. The final boundary conditions on $\theta$ and $\phi$ are subsequently varied in $1000$ equal increments to the desired values of $16.027$ deg and $1.1809$ deg respectively. The solution from a given iteration is used as the initial guess for the next iteration.
    \item In the second homotopy step, the final boundary condition on $v$ is fixed. For the first iteration in this step, this boundary condition is set to the final value of $v$ from the solution of the last iteration in step 1. This solution from step 1 also serves as the initial guess for the first iteration in step 2. The final boundary condition on $v$ is varied over $500$ iterations to the desired value of $0.1$ m/s.
    \item In the third homotopy step, $C_L$ for Segment 1 is varied in $500$ iterations from the current value of $0$ to the desired value of $0.25$ in equal increments. The initial guess for the first iteration in this step is the solution of the last iteration from step 2.
    \item In the fourth homotopy step, the initial boundary condition on $\psi$ is fixed. For the first iteration in this step, this boundary condition is set to the initial value of $\psi$ from the solution of the last iteration in step 3. This solution from step 3 also serves as the initial guess for the first iteration in step 4. The initial boundary condition on $\psi$ is varied over $500$ iterations to the desired value of $90$ degrees.
    \item In the fifth and final homotopy step, the slope parameters $s$ and $\zeta$ of the sigmoid and the hyperbolic tangent functions are varied over $1000$ iterations to $40,000$ from their current values of $100$ and $1$, respectively. The solution of the last iteration from step 4 serves as the initial guess for the first iteration of step 5. The solution of the last iteration of step 5 is the solution of the desired TPBVP.
\end{enumerate}

\begin{figure}[!h]
    \centering
    {\includegraphics[width=4in]{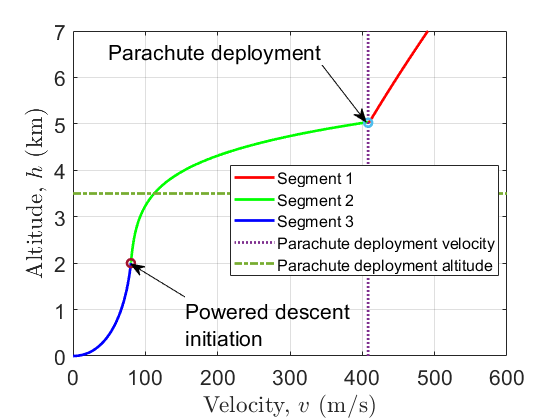}}
    \caption{\label{fig:hvZoomed_GRASHS_Profile1}
    Altitude vs. atmospheric-relative velocity for GRASHS solution for mission profile 1, zoomed in.}
\end{figure}

Figure \ref{fig:hvZoomed_GRASHS_Profile1} illustrates the variation of altitude $h$ as a function of atmospheric-relative velocity $v$ in the GRASHS solution. It is clear from this plot that the vehicle decelerated to the parachute deployment velocity ($v_P$) of $408$ m/s before descending to the parachute deployment altitude ($h_P$) of $3.5$ km, as indicated by the corner point at the vertical dash-line depicting $v_P$.

It is important to note that the GRASHS solution process neither entails assumptions not possesses knowledge about whether the vehicle first decelerates to $v_P$ or descends to $h_P$. This information is in fact already embedded into the equations of motion in Eq. \eqref{eqn:eom_GRASHS} through the sigmoid and hyperbolic tangent functions. The activation of appropriate flight segments, the specific trigger events (such as whether $v<v_P$ or $h<h_P$ activates parachute deployment), and the enforcement of the pertinent interior-point boundary conditions at the flight segment transition points are now merely a consequence of Eq. \eqref{eqn:eom_GRASHS} and are automatically and implicitly handled.

With the knowledge provided by the GRASHS solution that the parachute deployment was triggered by $v < v_P$, the trajectory is solved again using the RASHS approach by removing $h_P$ from the conditions listed in Table \ref{tab:flightSegmentConditions}, thereby eliminating the OR logic in Segment 2. Therefore, the DNF in Eq. \eqref{eqn:DNF_GRASHS} collapses to:

\begin{equation}
    \label{eqn:DNF_RASHS}
    \begin{gathered}
        \textrm{Segment 1: } \left(\left(\frac{v_P}{v\left(0\right)} - V\right) < 0\right) \\
        \textrm{Segment 2: } \left(\left(V - \frac{v_P}{v\left(0\right)}\right) < 0\right) \cdot \left(\left(\frac{h_{PDI}}{h\left(0\right)} - \hbar\right) < 0 \right) \\
        \textrm{Segment 3: } \left(\left(\hbar - \frac{h_{PDI}}{h\left(0\right)}\right) < 0\right)
    \end{gathered}
\end{equation}

\newpage Applying the RASHS approach, the equations of motion and the Lagrangian are as follows:

\begin{equation}
    \label{eqn:eom_RASHS}
    \textbf{f} = \left[\left(\frac{1}{1 + e^{s \left(\frac{v_P}{v\left(0\right)} - V\right)}}\right)\right] \textbf{f}_1 +
    \left[\left(\frac{1}{1 + e^{s \left(V - \frac{v_P}{v\left(0\right)}\right)}}\right)
    \left(\frac{1}{1 + e^{s \left(\frac{h_{PDI}}{h\left(0\right)} - \hbar\right)}}\right)\right] \textbf{f}_2
    +
    \left[\left(\frac{1}{1 + e^{s \left(\hbar - \frac{h_{PDI}}{h\left(0\right)} \right)}}\right)\right] \textbf{f}_3
\end{equation}

\begin{equation}
    \label{eqn:lagrangian_RASHS}
    \mathcal{L} = \left[\left(\frac{1}{1 + e^{s \left(\frac{v_P}{v\left(0\right)} - V\right)}}\right)\right] \mathcal{L}_1 +
    \left[\left(\frac{1}{1 + e^{s \left(V - \frac{v_P}{v\left(0\right)}\right)}}\right)
    \left(\frac{1}{1 + e^{s \left(\frac{h_{PDI}}{h\left(0\right)} - \hbar\right)}}\right)\right] \mathcal{L}_2
    +
    \left[\left(\frac{1}{1 + e^{s \left(\hbar - \frac{h_{PDI}}{h\left(0\right)} \right)}}\right)\right] \mathcal{L}_3
\end{equation}

The trajectory for the RASHS approach is solved using the same homotopy steps and initial guess generation technique (that employs \textit{ode45}) used for GRASHS, with the only exception that homotopy on $\zeta$ does not apply for RASHS because of the absence of the hyperbolic tangent function.

For comparison, the MPBVP from Eq. \eqref{eqn:mpbvp_necessary_conditions} is also solved by implementing $v_P$ and ignoring $h_P$ in the interior-point boundary conditions using the apriori knowledge gained from the GRASHS solution that the vehicle first decelerates to $v_P$ before descending to $h_P$. MATLAB \textit{bvp4c} is employed again to solve the MPBVP using the GRASHS solution as the initial guess.

\begin{figure}[!]
	\centering
	\begin{subfigure}[t]{0.49\textwidth}
		\centering
		\includegraphics[width=\textwidth]{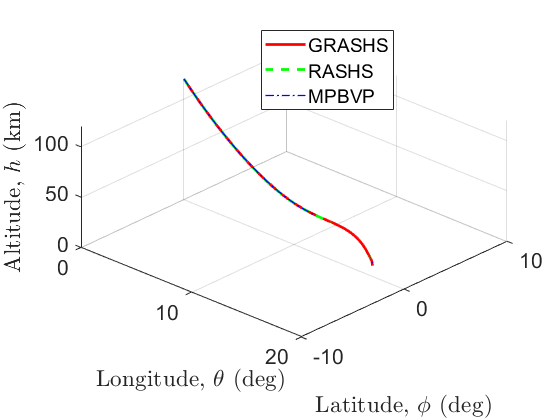}
		\caption{\label{fig:trajectory_Profile1}Physical trajectory.}
	\end{subfigure}
	\begin{subfigure}[t]{0.49\textwidth}
		\centering
		\includegraphics[width=\textwidth]{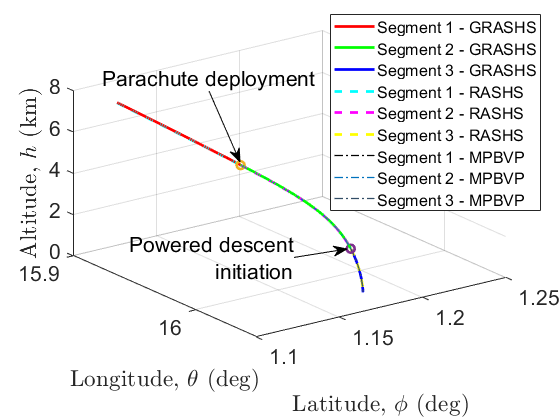}
		\caption{\label{fig:trajectoryZoomed_Profile1}Physical trajectory zoomed in.}
	\end{subfigure}
    \begin{subfigure}[t]{0.49\textwidth}
		\centering
		\includegraphics[width=\textwidth]{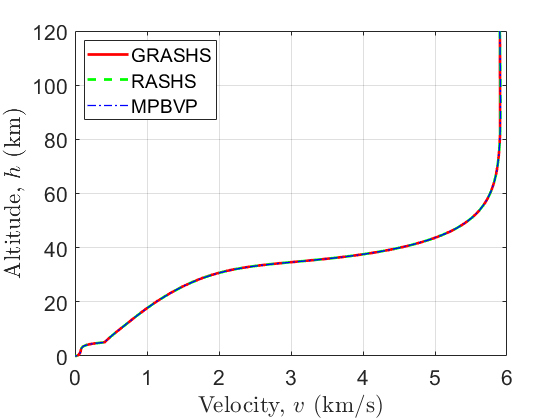}
		\caption{\label{fig:hv_Profile1}Altitude vs. atmospheric-relative velocity}
	\end{subfigure}
    \begin{subfigure}[t]{0.49\textwidth}
		\centering
		\includegraphics[width=\textwidth]{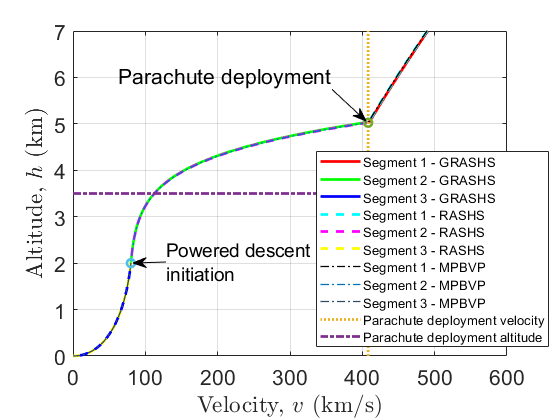}
		\caption{\label{fig:hvZoomed_Profile1}Altitude vs. atmospheric-relative velocity, zoomed in}
	\end{subfigure}
    \begin{subfigure}[t]{0.49\textwidth}
		\centering
		\includegraphics[width=\textwidth]{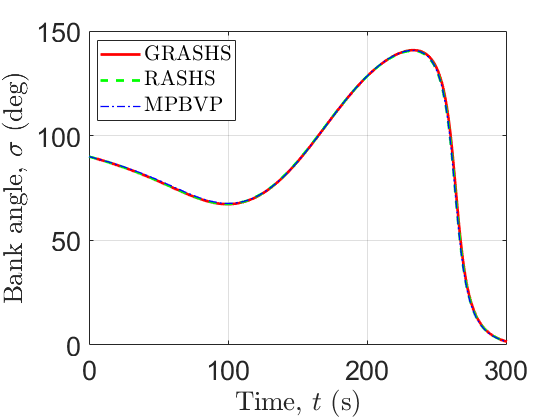}
		\caption{\label{fig:bankAngle_Profile1}Bank angle vs. time}
	\end{subfigure}
    \begin{subfigure}[t]{0.49\textwidth}
		\centering
		\includegraphics[width=\textwidth]{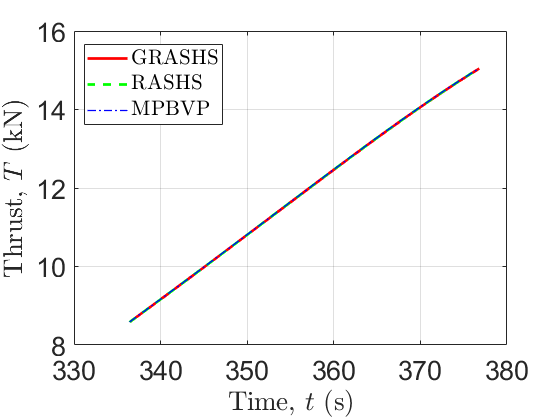}
		\caption{\label{fig:thrust_Profile1}Thrust vs. time}
	\end{subfigure}
	\caption{\label{fig:StatesAndControl_Profile1}Comparison of trajectory and control for mission profile 1.}
\end{figure}

\begin{figure}[!h]
    \centering
    {\includegraphics[width=3.5in]{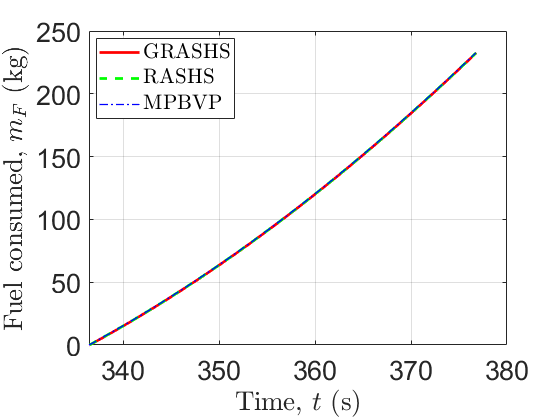}}
    \caption{\label{fig:fuel_Profile1}Propellant consumption vs. time for mission profile 1.}
\end{figure}

\begin{figure}[!]
	\centering
	\begin{subfigure}[t]{0.37\textwidth}
		\centering
		\includegraphics[width=\textwidth]{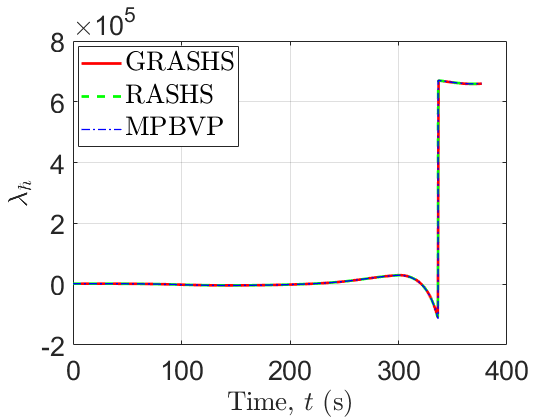}
		\caption{\label{fig:lam_hbar_Profile1}$\lambda_{\hbar}$ vs. $t$.}
	\end{subfigure}
	\begin{subfigure}[t]{0.37\textwidth}
		\centering
		\includegraphics[width=\textwidth]{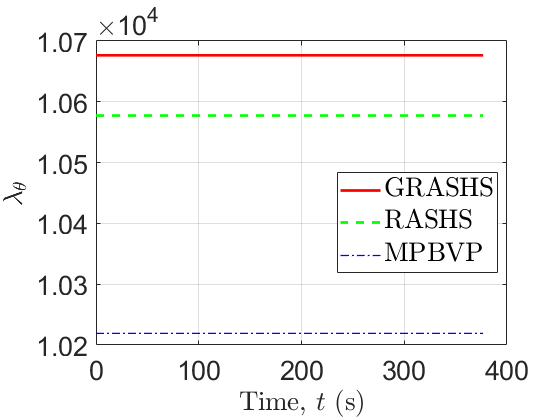}
		\caption{\label{fig:lam_theta_Profile1}$\lambda_{\theta}$ vs. $t$.}
	\end{subfigure}
    \begin{subfigure}[t]{0.37\textwidth}
		\centering
		\includegraphics[width=\textwidth]{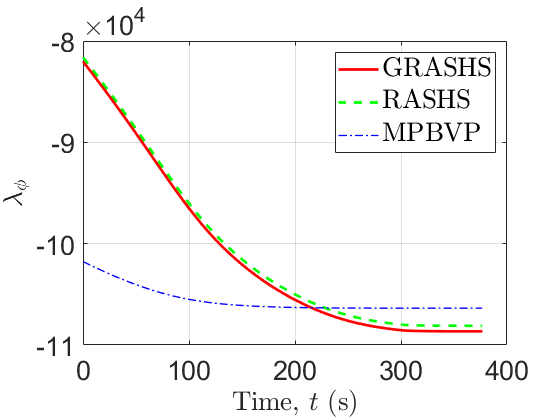}
		\caption{\label{fig:lam_phi_Profile1}$\lambda_{\phi}$ vs. $t$.}
	\end{subfigure}
    \begin{subfigure}[t]{0.37\textwidth}
		\centering
		\includegraphics[width=\textwidth]{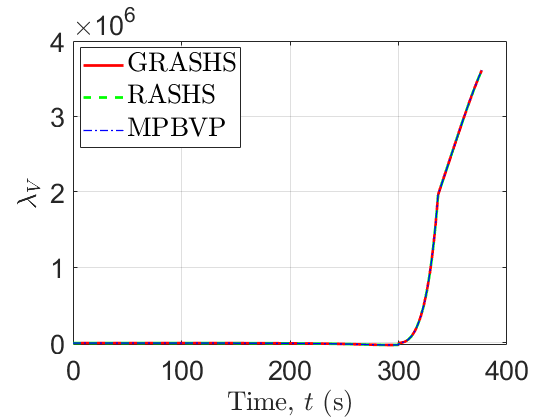}
		\caption{\label{fig:lam_V_Profile1}$\lambda_{V}$ vs. $t$.}
	\end{subfigure}
    \begin{subfigure}[t]{0.37\textwidth}
		\centering
		\includegraphics[width=\textwidth]{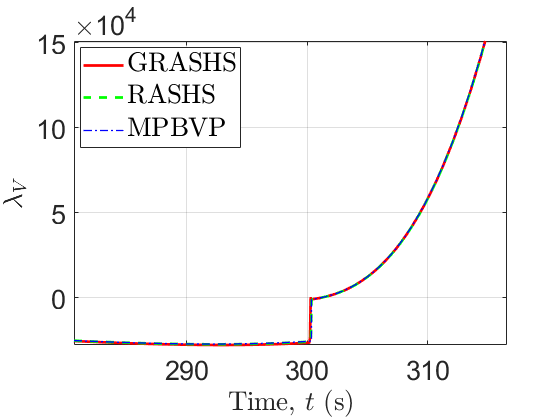}
		\caption{\label{fig:lam_V_Zoomed_Profile1}$\lambda_{V}$ vs. $t$, zoomed in.}
	\end{subfigure}
    \begin{subfigure}[t]{0.37\textwidth}
		\centering
		\includegraphics[width=\textwidth]{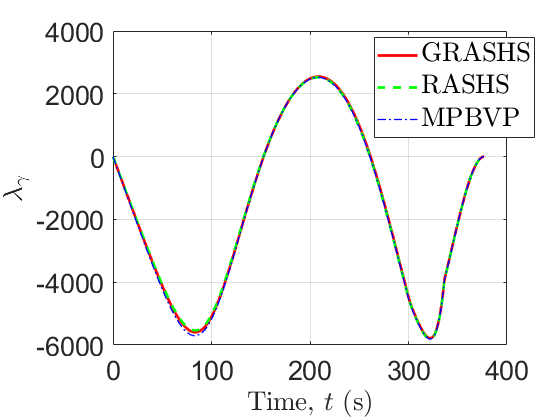}
		\caption{\label{fig:lam_gam_profile1}$\lambda_{\gamma}$ vs. $t$.}
	\end{subfigure}
    \begin{subfigure}[t]{0.37\textwidth}
		\centering
		\includegraphics[width=\textwidth]{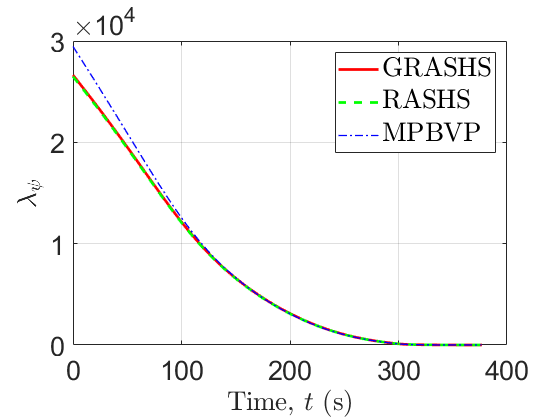}
		\caption{\label{fig:lam_psi_Profile1}$\lambda_{\psi}$ vs. $t$.}
	\end{subfigure}
    \begin{subfigure}[t]{0.37\textwidth}
		\centering
		\includegraphics[width=\textwidth]{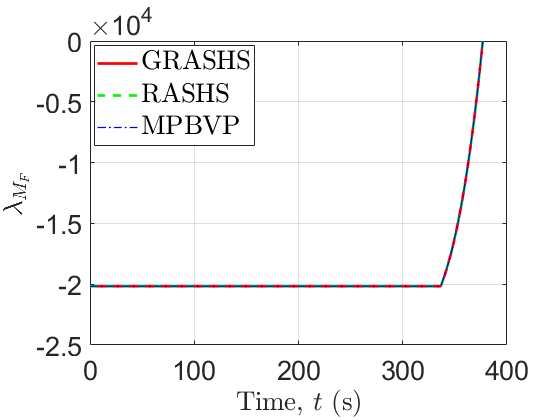}
		\caption{\label{fig:lam_M_F_Profile1}$\lambda_{M_F}$ vs. $t$.}
	\end{subfigure}
	\caption{\label{fig:costates_Profile1}Comparison of co-state history for mission profile 1.}
\end{figure}

Figure \ref{fig:trajectory_Profile1} compares the physical three-dimensional trajectories from GRASHS, RASHS and MPBVP. Figure \ref{fig:trajectoryZoomed_Profile1} shows the same plot zoomed in to illustrate the parachute deployment and the PDI events. The plots demonstrate that the physical trajectories generated by GRASHS, RASHS and MPBVP are consistent. This is corroborated by the parachute deployment, PDI and touchdown times, as summarized in Table \ref{tab:eventTimes_Profile1}.

\begin{table}[!h]
    \centering
    \caption{\label{tab:eventTimes_Profile1}Parachute deployment, PDI and touchdown times for mission profile 1 derived from GRASHS, RASHS and MPBVP solutions.}
    \begin{tabular}{lccc}
        \hline
        Method & Parachute deployment & PDI & Touchdown \\\hline
        GRASHS & 300.3052 \textrm{s} & 336.4748 \textrm{s} & 376.7966 \textrm{s} \\
        RASHS & 300.3600 \textrm{s} & 336.4492 \textrm{s} & 376.7730 \textrm{s} \\
        MPBVP & 300.3537 \textrm{s} & 336.4169 \textrm{s} & 376.7456 \textrm{s} \\
        \hline
    \end{tabular}
\end{table}

Figure \ref{fig:hv_Profile1} compares the plots of altitude as a function of atmospheric-relative velocity from the three solutions, and Fig. \ref{fig:hvZoomed_Profile1} shows the same plot zoomed in. As before, the solutions show good consistency. The parachute deployment event is indicated by the corner point in Fig. \ref{fig:hvZoomed_Profile1}, at the vertical dash-line depicting $v_P$. The GRASHS and RASHS solutions are actually smooth in this section of the plot. However, because of the steep slope parameters $s$ and $\zeta$, the smooth transitions appear to be a corner point. This renders the errors associated with the smoothing from sigmoid and hyperbolic tangent functions to be negligible for practical purposes. Although not apparent in the plot, there is also a corner point at PDI because of the change in deceleration resulting from change in mass, $C_D$ and the application of thrust.

The bank angle (Fig. \ref{fig:bankAngle_Profile1}) and thrust (Fig. \ref{fig:thrust_Profile1}) histories are also consistent across the three solutions. As expected, the vehicle banks to the left, as indicated by the positive value, to turn northbound because the vehicle is flying due East and the touchdown location is towards northeast. Consistent with the thrust, Fig. \ref{fig:fuel_Profile1} illustrates the propellant consumption over time, and shows that the three solutions match.

Figure \ref{fig:costates_Profile1} compares the co-state history from the three solutions. The histories of $\lambda_{\hbar}$ (Fig. \ref{fig:lam_hbar_Profile1}), $\lambda_V$ (Fig. \ref{fig:lam_V_Profile1} and \ref{fig:lam_V_Zoomed_Profile1}), $\lambda_{\gamma}$ (Fig. \ref{fig:lam_gam_profile1}) and $\lambda_{M_F}$ (Fig. \ref{fig:lam_M_F_Profile1}) are consistent across the three solutions. Of particular interest are the histories of $\lambda_{\hbar}$ and $\lambda_V$, which are the co-state histories of nondimensional altitude and atmospheric-relative velocity. Because parachute deployment is velocity-triggered, the implication is that the velocity is fixed at parachute deployment ($v\left(t_P\right) = v_P$). Therefore, as predicted by Eq. \eqref{eqn:mpbvp_necessary_conditions} (specifically, the seventh sub-equation that states ${\bm{\lambda}^T}\left(t^-_{k\neq m}\right) = {\bm{\lambda}^T}\left(t^+_{k\neq m}\right) + \bm{\Pi}_k^T \frac{\partial \bm{\Psi}_k}{\partial \textbf{X}\left(t_k\right)}$), $\lambda_V$ is expected to jump at parachute deployment. This jump can be observed in Fig. \ref{fig:lam_V_Zoomed_Profile1}, which is essentially a zoomed-in version of Fig. \ref{fig:lam_V_Profile1}. The jump in $\lambda_V$ is the same in all three solutions.

Similarly, because the altitude is fixed at PDI ($h\left(t_{PDI}\right) = h_{PDI}$), $\lambda_{\hbar}$ is expected to jump at PDI. This can be observed in Fig. \ref{fig:lam_hbar_Profile1}, and the jump from all three solutions match one another. It must be noted that for GRASHS (and RASHS), these jumps are automatically governed by the equations of motion and Lagrangian in Eqs. \eqref{eqn:eom_GRASHS} and \eqref{eqn:lagrangian_GRASHS} (Eqs. \eqref{eqn:eom_RASHS} and \eqref{eqn:lagrangian_RASHS} for RASHS) because they embed the conditions that govern the active flight segment (Table \ref{tab:flightSegmentConditions} for GRASHS and Table \ref{tab:flightSegmentConditions} with $h_P$ ignored for RASHS). Therefore, these jumps do not have to be explicitly calculated. Because of the sigmoid and the hyperbolic tangent functions, these jumps are actually smooth transitions in these solutions, but appear discrete because of the large values of $s$ and $\zeta$. As a result, the errors in the states are minimal and the GRASHS and RASHS solutions are consistent with that of the MPBVP for practical purposes, as evidenced by Fig. \ref{fig:StatesAndControl_Profile1} and \ref{fig:fuel_Profile1}. In the MPBVP, the jumps in the co-states must be explicitly calculated as part of the solution process, which would have been a challenging task if a good initial guess (such as the GRASHS solution) was not available.

The co-states $\lambda_{\theta}$ (Fig. \ref{fig:lam_theta_Profile1}), $\lambda_{\phi}$ (Fig. \ref{fig:lam_phi_Profile1}) and $\lambda_{\psi}$ (Fig. \ref{fig:lam_psi_Profile1}) exhibit some deviation. As explained in \cite{rashs}, it may be inferred that the errors arising from the smoothing operation effected by the sigmoid and hyperbolic tangent functions manifest in these co-states. These deviations were observed to diminish with increase in the values of $s$ and $\zeta$. However, it must be noted that these are the optimal solutions for the smoothed problems in GRASHS and RASHS because the necessary conditions in Eq. \eqref{eqn:tpbvp_necessary_conditions} are satisfied. Despite these deviations, the states are practically consistent with those of the MPBVP.

The results in this section demonstrated that the GRASHS approach was able to solve the multi-phase EDL trajectory with no apriori knowledge about whether the parachute deployment event was velocity-triggered or altitude-triggered. The fact that this solution was consistent with those of RASHS and MPBVP, which relied on apriori assumption about the trigger event, demonstrated the advantage of the GRASHS approach over the state-of-the-art for indirect methods. If the apriori assumptions were to change, the RASHS and MPBVP solutions would become invalid because Eqs. \eqref{eqn:DNF_RASHS}, \eqref{eqn:eom_RASHS} and \eqref{eqn:lagrangian_RASHS} would no longer apply, as will be demonstrated in Sec. \ref{sec:missionProfile2}.

\subsection{Mission Profile 2: High Parachute Deployment Altitude}
\label{sec:missionProfile2}
In this section, the optimal trajectory is calculated again for the same Mars EDL scenario, but with the parachute deployment altitude $h_P$ increased to $6.5$ km. For the GRASHS approach, the trajectory is solved using exactly the same initial guess generation technique and homotopy steps as described in Sec. \ref{sec:missionProfile1}. It is important to note that because GRASHS does not entail any assumptions about the trigger events, the same equations of motion and cost functional from Eqs. \eqref{eqn:eom_GRASHS} and \eqref{eqn:lagrangian_GRASHS} in Sec. \ref{sec:missionProfile1} can be reused, with the updated $h_P$ value.

\begin{figure}[!h]
    \centering
    {\includegraphics[width=4in]{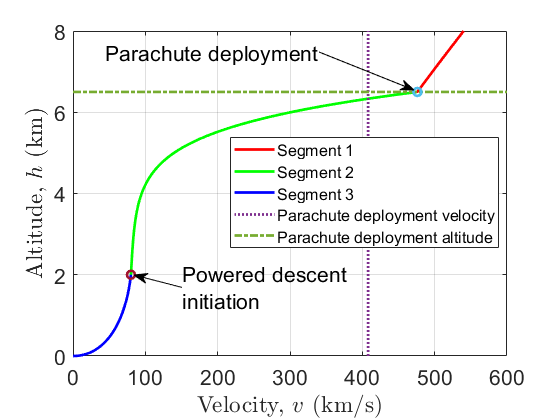}}
    \caption{\label{fig:hvZoomed_GRASHS_Profile2}
    Altitude vs. atmospheric-relative velocity for GRASHS solution for mission profile 2, zoomed in.}
\end{figure}

This time, the altitude vs. atmospheric-relative velocity plot, illustrated in Fig. \ref{fig:hvZoomed_GRASHS_Profile2}, shows that the vehicle descends to the parachute deployment altitude ($h_P$) of $6.5$ km before decelerating to the parachute deployment velocity ($v_P$) of $408$ m/s, as indicated by the corner point at the horizontal dash-dot-line depicting $h_P$. Again, it must be noted that this information is not available to the GRASH solution process, and it not even required to be known in the GRASHS framework. Instead, this information is already embedded into the equations of motion in Eq. \eqref{eqn:eom_GRASHS} through the sigmoid and hyperbolic tangent functions and the parachute descent segment is automatically activated based on which trigger event is hit first. In this particular scenario, it so happens that the equations of motion dictate that the vehicle hits $h_P$ first.

Having gained knowledge from the GRASHS solution that the vehicle hits $h_P$ first, the trajectory is solved again using RASHS by removing $v_P$ from the conditions listed in Table \ref{tab:flightSegmentConditions}, thereby eliminating the OR logic in Segment 2. Accordingly, this time, the DNF for RASHS collapses to:

\begin{equation}
    \label{eqn:DNF_RASHS_Profile2}
    \begin{gathered}
        \textrm{Segment 1: } \left(\left(\frac{h_P}{h\left(0\right)} - \hbar\right) < 0\right) \\
        \textrm{Segment 2: } \left(\left(\hbar - \frac{h_P}{h\left(0\right)}\right) < 0\right) \cdot \left(\left(\frac{h_{PDI}}{h\left(0\right)} - \hbar\right) < 0 \right) \\
        \textrm{Segment 3: } \left(\left(\hbar - \frac{h_{PDI}}{h\left(0\right)}\right) < 0\right)
    \end{gathered}
\end{equation}

Consequently, the equations of motion and the Lagrangian for the RASHS approach are as follows:

\begin{equation}
    \label{eqn:eom_RASHS_Profile2}
    \textbf{f} = \left[\left(\frac{1}{1 + e^{s \left(\frac{h_P}{h\left(0\right)} - \hbar\right)}}\right)\right] \textbf{f}_1 +
    \left[\left(\frac{1}{1 + e^{s \left(\hbar - \frac{h_P}{h\left(0\right)}\right)}}\right)
    \left(\frac{1}{1 + e^{s \left(\frac{h_{PDI}}{h\left(0\right)} - \hbar\right)}}\right)\right] \textbf{f}_2
    +
    \left[\left(\frac{1}{1 + e^{s \left(\hbar - \frac{h_{PDI}}{h\left(0\right)} \right)}}\right)\right] \textbf{f}_3
\end{equation}

\begin{multline}
    \label{eqn:lagrangian_RASHS_Profile2}
    \mathcal{L} = \left[\left(\frac{1}{1 + e^{s \left(\frac{h_P}{h\left(0\right)} - \hbar\right)}}\right)\right] \mathcal{L}_1 +
    \left[\left(\frac{1}{1 + e^{s \left(\hbar - \frac{h_P}{h\left(0\right)}\right)}}\right)
    \left(\frac{1}{1 + e^{s \left(\frac{h_{PDI}}{h\left(0\right)} - \hbar\right)}}\right)\right] \mathcal{L}_2
    + \\
    \left[\left(\frac{1}{1 + e^{s \left(\hbar - \frac{h_{PDI}}{h\left(0\right)} \right)}}\right)\right] \mathcal{L}_3
\end{multline}

Note that Eqs. \eqref{eqn:eom_RASHS_Profile2} and \eqref{eqn:lagrangian_RASHS_Profile2} are different from Eqs. \eqref{eqn:eom_RASHS} and \eqref{eqn:lagrangian_RASHS} in Sec. \ref{sec:missionProfile1}. This shows that the RASHS approach requires reformulation of the equations of motion and the Lagrangian if the underlying assumptions about the trigger events change.

The TPBVP for the RASHS approach is also solved using the same initial guess generation technique and homotopy steps used in the GRASHS approach in Sec. \ref{sec:missionProfile1}. Finally, the MPBVP in Eq. \eqref{eqn:mpbvp_necessary_conditions} is solved using the new GRASHS solution as the initial guess. This time, the MPBVP implements $h_P$ and ignores $v_P$ in the interior-point boundary conditions using the apriori knowledge gained from the new GRASHS solution that the vehicle descends to $h_P$ before decelerating to $v_P$.

\begin{figure}[!]
	\centering
	\begin{subfigure}[t]{0.49\textwidth}
		\centering
		\includegraphics[width=\textwidth]{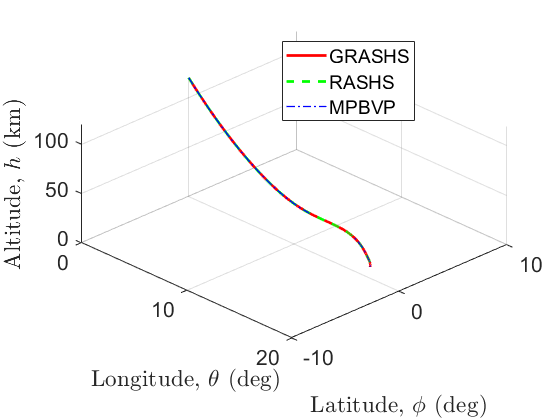}
		\caption{\label{fig:trajectory_Profile2}Physical trajectory.}
	\end{subfigure}
	\begin{subfigure}[t]{0.49\textwidth}
		\centering
		\includegraphics[width=\textwidth]{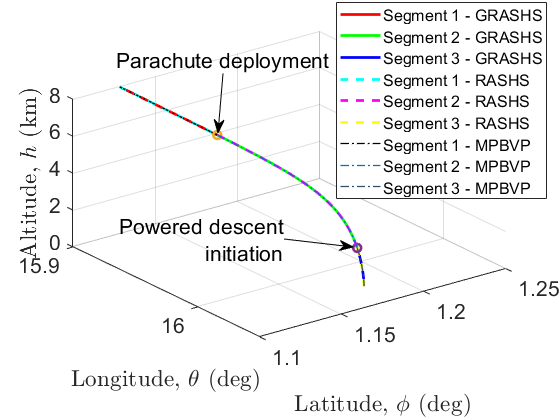}
		\caption{\label{fig:trajectoryZoomed_Profile2}Physical trajectory zoomed in.}
	\end{subfigure}
    \begin{subfigure}[t]{0.49\textwidth}
		\centering
		\includegraphics[width=\textwidth]{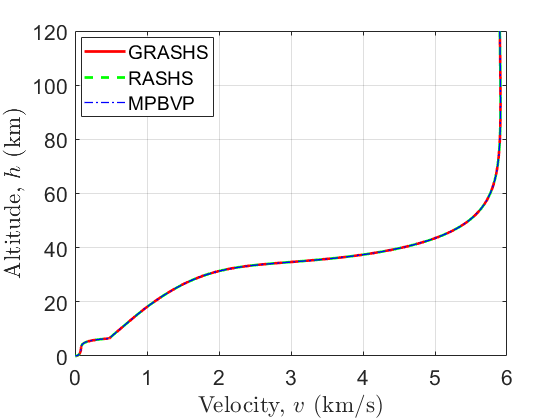}
		\caption{\label{fig:hv_Profile2}Altitude vs. atmospheric-relative velocity}
	\end{subfigure}
    \begin{subfigure}[t]{0.49\textwidth}
		\centering
		\includegraphics[width=\textwidth]{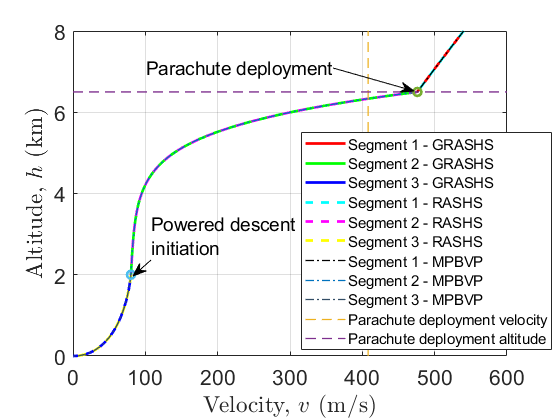}
		\caption{\label{fig:hvZoomed_Profile2}Altitude vs. atmospheric-relative velocity, zoomed in}
	\end{subfigure}
    \begin{subfigure}[t]{0.49\textwidth}
		\centering
		\includegraphics[width=\textwidth]{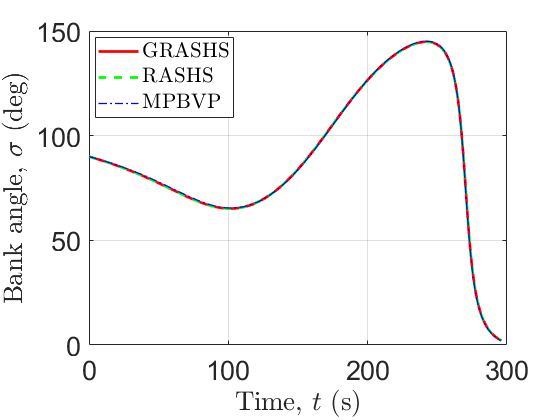}
		\caption{\label{fig:bankAngle_Profile2}Bank angle vs. time}
	\end{subfigure}
    \begin{subfigure}[t]{0.49\textwidth}
		\centering
		\includegraphics[width=\textwidth]{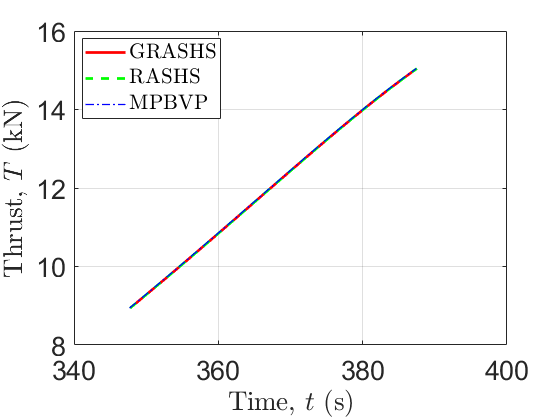}
		\caption{\label{fig:thrust_Profile2}Thrust vs. time}
	\end{subfigure}
	\caption{\label{fig:StatesAndControl_Profile2}Comparison of trajectory and control for mission profile 2.}
\end{figure}

\begin{figure}[!h]
    \centering
    {\includegraphics[width=3.5in]{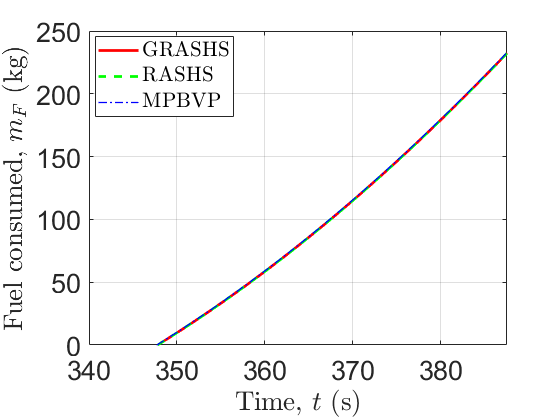}}
    \caption{\label{fig:fuel_Profile2}Propellant consumption vs. time for mission profile 2.}
\end{figure}

Figures \ref{fig:StatesAndControl_Profile2} and \ref{fig:fuel_Profile2} compare the solutions of GRASHS, RASHS and MPBVP. Figures \ref{fig:trajectory_Profile2} and \ref{fig:trajectoryZoomed_Profile2} show the physical three-dimensional trajectories, with the parachute deployment and PDI events illustrated in the latter. The physical trajectories closely match one another, which is further corroborated by the parachute deployment, PDI and touchdown times summarized in Table \ref{tab:eventTimes_Profile2}.

\begin{table}[!h]
    \centering
    \caption{\label{tab:eventTimes_Profile2}Parachute deployment, PDI and touchdown times for mission profile 2.}
    \begin{tabular}{lccc}
        \hline
        Event & Parachute deployment & PDI & Touchdown \\\hline
        GRASHS & 295.9386 \textrm{s} & 347.7960 \textrm{s} & 387.5228 \textrm{s} \\
        RASHS & 295.9402 \textrm{s} & 347.7949 \textrm{s} & 387.5229 \textrm{s} \\
        MPBVP & 295.8594 \textrm{s} & 347.7101 \textrm{s} & 387.4404 \textrm{s} \\
        \hline
    \end{tabular}
\end{table}

Figures \ref{fig:hv_Profile2} and \ref{fig:hvZoomed_Profile2} show the altitude as a function of atmospheric-relative velocity. As before, the solutions are consistent. The parachute deployment event is indicated by the corner point in Fig. \ref{fig:hvZoomed_Profile2} at the horizontal dash-dot-line representing $h_P$. Also, the bank angle (Fig. \ref{fig:bankAngle_Profile2}), thrust (Fig. \ref{fig:thrust_Profile2}) and propellant consumption (\ref{fig:fuel_Profile2}) are consistent across the three solutions and the vehicle banks to the left.

\begin{figure}[!]
	\centering
	\begin{subfigure}[t]{0.37\textwidth}
		\centering
		\includegraphics[width=\textwidth]{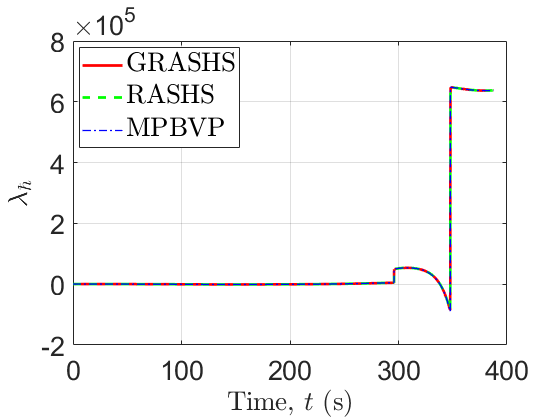}
		\caption{\label{fig:lam_hbar_Profile2}$\lambda_{\hbar}$ vs. $t$.}
	\end{subfigure}
	\begin{subfigure}[t]{0.37\textwidth}
		\centering
		\includegraphics[width=\textwidth]{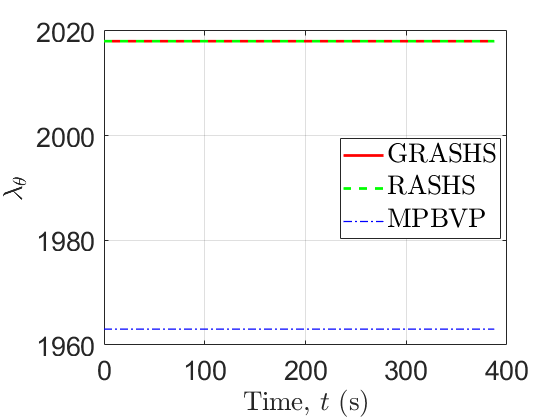}
		\caption{\label{fig:lam_theta_Profile2}$\lambda_{\theta}$ vs. $t$.}
	\end{subfigure}
    \begin{subfigure}[t]{0.37\textwidth}
		\centering
		\includegraphics[width=\textwidth]{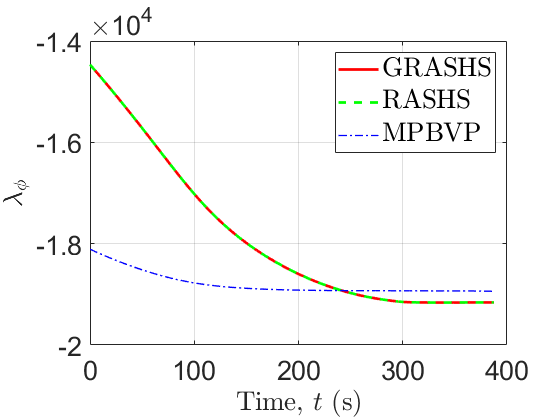}
		\caption{\label{fig:lam_phi_Profile2}$\lambda_{\phi}$ vs. $t$.}
	\end{subfigure}
    \begin{subfigure}[t]{0.37\textwidth}
		\centering
		\includegraphics[width=\textwidth]{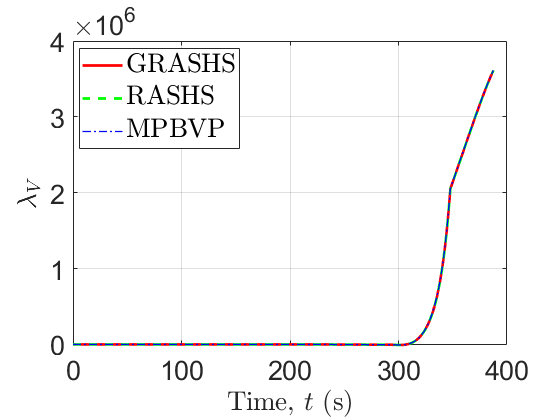}
		\caption{\label{fig:lam_V_Profile2}$\lambda_{V}$ vs. $t$.}
	\end{subfigure}
    \begin{subfigure}[t]{0.37\textwidth}
		\centering
		\includegraphics[width=\textwidth]{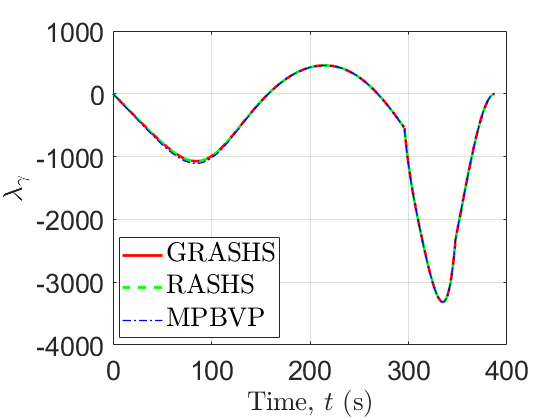}
		\caption{\label{fig:lam_gam_profile2}$\lambda_{\gamma}$ vs. $t$.}
	\end{subfigure}
    \begin{subfigure}[t]{0.37\textwidth}
		\centering
		\includegraphics[width=\textwidth]{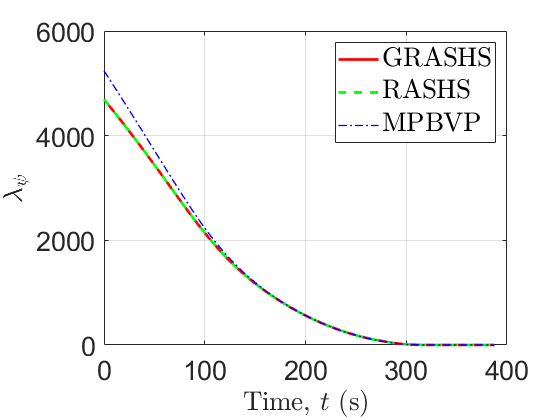}
		\caption{\label{fig:lam_psi_Profile2}$\lambda_{\psi}$ vs. $t$.}
	\end{subfigure}
    \begin{subfigure}[t]{0.37\textwidth}
		\centering
		\includegraphics[width=\textwidth]{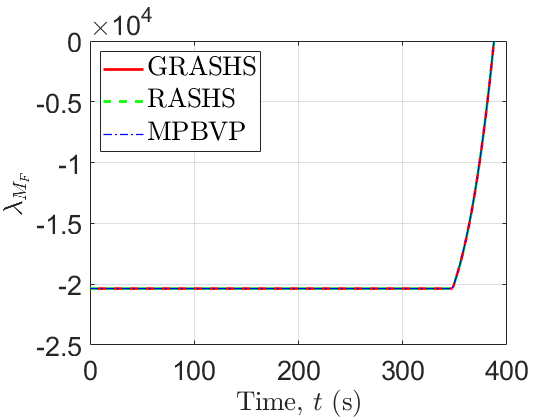}
		\caption{\label{fig:lam_M_F_Profile2}$\lambda_{M_F}$ vs. $t$.}
	\end{subfigure}
	\caption{\label{fig:costates_Profile2}Comparison of co-state history for mission profile 2.}
\end{figure}

Figure \ref{fig:costates_Profile2} shows the comparison of the co-state histories. As before, $\lambda_{\hbar}$ (Fig. \ref{fig:lam_hbar_Profile2}), $\lambda_V$ (Fig. \ref{fig:lam_V_Profile2}), $\lambda_{\gamma}$ (Fig. \ref{fig:lam_gam_profile2}) and $\lambda_{M_F}$ (Fig. \ref{fig:lam_M_F_Profile2}) match well. However, in this mission profile, $\lambda_{\hbar}$ exhibits a jump twice (Fig. \ref{fig:lam_hbar_Profile2}) because the altitude is fixed at both parachute deployment ($h = h_P$) and PDI ($h = h_{PDI}$) events. Also, unlike before, $\lambda_V$ does not exhibit any jump (Fig. \ref{fig:lam_V_Profile2}) because neither parachute deployment nor PDI is associated with velocity. This is because the parachute deployed before hitting $v_P$, and hence, $v = v_P$ no longer constitutes an interior-point boundary condition at parachute deployment event. Therefore, the value of $v$ is free at both parachute deployment and PDI. As a result, ${\bm{\lambda}^T}\left(t^-_{k\neq m}\right) = {\bm{\lambda}^T}\left(t^+_{k\neq m}\right) + \bm{\Pi}_k^T \frac{\partial \bm{\Psi}_k}{\partial \textbf{X}\left(t_k\right)}$ in Eq. \eqref{eqn:mpbvp_necessary_conditions} collapses to ${\bm{\lambda}^T}\left(t^-_{k\neq m}\right) = {\bm{\lambda}^T}\left(t^+_{k\neq m}\right)$ for $\lambda_V$ at both $t_P$ and $t_{PDI}$, rendering it continuous for the entire trajectory.

Finally, $\lambda_{\theta}$ (Fig. \ref{fig:lam_theta_Profile2}), $\lambda_{\phi}$ (Fig. \ref{fig:lam_phi_Profile2}) and $\lambda_{\psi}$ (Fig. \ref{fig:lam_psi_Profile2}) exhibit deviations as seen in Fig. \ref{fig:costates_Profile1} in Sec. \ref{sec:missionProfile1}, which can be attributed to the errors introduced by the smoothing operations from sigmoid and hyperbolic tangent functions. Nevertheless, to reiterate, the GRASHS and RASHS solutions are in fact optimal because they satisfy the necessary conditions of optimality in Eq. $\ref{eqn:tpbvp_necessary_conditions}$. For practical purposes, the states calculated using the three methods match well, as evidenced by Fig. \ref{fig:StatesAndControl_Profile2} and \ref{fig:fuel_Profile2}.

This section demonstrated that although the trigger condition for parachute deployment changed from velocity to altitude, the GRASHS solution was still able to reuse the same equations of motion and Lagrangian defined in Eqs. \eqref{eqn:eom_GRASHS} and \eqref{eqn:lagrangian_GRASHS} from Sec. \ref{sec:missionProfile1}, and continued to make no assumptions about the trigger event. This is unlike the RASHS and MPBVP solutions, which changed their underlying assumptions from Sec. \ref{sec:missionProfile1} based on the new GRASHS solution from this section. Accordingly, for RASHS, the equations of motion and Lagrangian had to be updated accordingly, as evidenced by Eqs. \eqref{eqn:DNF_RASHS_Profile2}, \eqref{eqn:eom_RASHS_Profile2} and \eqref{eqn:lagrangian_RASHS_Profile2}. Additionally, the interior-point boundary conditions of the MPBVP also had to be updated based on the updated assumptions about the trigger event. Therefore, although both GRASHS and RASHS have advantages over the original MPBVP, GRASHS is a clear improvement over RASHS.

\section{Conclusion}
This investigation presented a methodology for the indirect multi-phase trajectory optimization framework to improve the RASHS approach to handle cases when the conditions governing the active flight segment are represented using an arbitrary discrete logic. Although the original RASHS formulation effeectively simplified the multi-phase trajectory optimization process by reducing the MPBVP representing the necessary conditions of optimality to an easier-to-handle TPBVP, it was designed to only handle AND logic. The new methodology, termed the Generalized Relaxed Autonomously Switched Hybrid System (GRASHS) approach, achieved the capability to handle the arbitrary discrete logic by first transforming the logic to DNF. The NOT operations in the DNF were handled by replacing the pertinent predicates with ones bearing the form $g_{i,j,k} < 0$. The minterms in the DNF were represented using products of sigmoid functions and the OR operation on the minterms was represented as the hyperbolic tangent of the summation of the minterms. Because the resultant equations of motion and Lagrangian of the path cost were continuous and differential for the entire trajectory, the necessary conditions of optimality constituted a TPBVP in a system of DAEs, while traditional approaches would have resulted in an MPBVP that is difficult to handle. Moreover, because every boolean expression can be represented solely using AND, OR and NOT logic, the GRASHS approach can handle any arbitrary discrete logic.

The GRASHS approach was demonstrated by applying it to a multi-phase Mars EDL trajectory optimization example, where the conditions governing the parachute descent segment were represented using a combination of AND and OR logic. Specifically, this segment was active when the velocity was below $v_P$ OR the altitude was below $h_P$, AND the altitude was above $h_{PDI}$. Two mission profiles were presented, where $h_P$ was set to $3.5$ km in the first profile, and $6.5$ km in the second profile. Consequently, parachute deployment was triggered by $v_P$ in the former, and $h_P$ in the latter. For comparison, the problem was also solved using the RASHS approach and the original MPBVP by eliminating the OR logic through assumptions about whether the vehicle decelerated to $v_P$ or descended to $h_P$ first. The states and control showed consistency across the three solutions for both profiles.

The results clearly demonstrated that both GRASHS and RASHS approaches simplified the design of the multi-phase EDL trajectory by reducing the MPBVP to a TPBVP. However, GRASHS was an improvement over the RASHS approach because the former did not entail any assumptions or possess knowledge about whether the vehicle first decelerated to $v_P$ or descended to $h_P$. This was evident from the fact that both mission profiles used the same equations of motion (Eq. \eqref{eqn:eom_GRASHS}) and Lagrangian (\eqref{eqn:lagrangian_GRASHS}). The knowledge of whether the vehicle hit $v_P$ or $h_P$ first was not required for the GRASHS approach because this information was already embedded into Eq. \eqref{eqn:eom_GRASHS}, which also automatically transitioned the flight to parachute descent upon encountering any one of the triggers. This was unlike the RASHS approach, which required this apriori knowledge. Based on this knowledge, the OR logic was eliminated by removing the $h_P$ terms in Sec. \ref{sec:missionProfile1} and $v_P$ terms in Sec. \ref{sec:missionProfile2}. Consequently, the equations of motion and Lagrangian were different for the two mission profiles (Eqs. \eqref{eqn:eom_RASHS}, \eqref{eqn:lagrangian_RASHS}, \eqref{eqn:eom_RASHS_Profile2},  and \eqref{eqn:lagrangian_RASHS_Profile2}), thereby effectively demonstrating that the RASHS approach required reformulation any time the underlying assumptions about the trigger events changed. The fact that the GRASHS solution was consistent with that of RASHS despite not making any assumptions about the trigger events highlighted its advantage over the latter.

Finally, as seen in the original RASHS approach, the jump in the co-states at the flight segment transition points in the GRASHS approach were also implicitly handled by the equations of motion and Lagrangian in Eqs. \eqref{eqn:eom_GRASHS} and \eqref{eqn:lagrangian_GRASHS} and did not entail any explicit calculation. This is because the GRASHS approach collapsed the MPBVP to a TPBVP by essentially embedding the interior-point boundary conditions into Eqs. \eqref{eqn:eom_GRASHS} and \eqref{eqn:lagrangian_GRASHS}.

\bibliographystyle{unsrtnat}
\bibliography{grashs}

\end{document}